\documentclass[12pt]{article}

\usepackage{amsmath}
\usepackage{amssymb}
\usepackage{dsfont} 

\usepackage{graphicx}

\usepackage[dvipsnames]{xcolor}

\renewcommand{\theequation}{\arabic{section}.\arabic{equation}}

\newcommand{\be}{\begin{equation}}
\newcommand{\bel}[1]{\begin{equation}\label{#1}}
\newcommand{\ee}{\end{equation}}
\newcommand{\bea}{\begin{eqnarray}}
\newcommand{\eea}{\end{eqnarray}}
\newcommand{\balign}{\begin{align}}
\newcommand{\ealign}{\end{align}}
\newcommand{\ba}{\begin{array}}
\newcommand{\ea}{\end{array}}
\newcommand{\bfig}{\begin{figure}}
\newcommand{\efig}{\end{figure}}

\newcommand{\eref}[1]{(\ref{#1})}

\allowdisplaybreaks

\newcommand{\bra}[1]{\mbox{$\langle \, {#1}\, |$}}
\newcommand{\ket}[1]{\mbox{$| \, {#1}\, \rangle$}}
\newcommand{\exval}[1]{\mbox{$\langle \, {#1}\, \rangle$}}
\newcommand{\inprod}[2]{\mbox{$\langle \, {#1} \, | \, {#2} \, \rangle$}}

\newcommand{\E}{{\mathbf E}}

\newcommand{\bfx}{\mathbf{x}}
\newcommand{\bfy}{\mathbf{y}}

\newcommand{\rme}{\mathrm{e}}

\newcommand{\half}{\frac{1}{2}}

\newcommand{\R}{{\mathbb R}}

\newcommand{\Z}{{\mathbb Z}}
\newcommand{\N}{{\mathbb N}}

\newcommand{\bfeta}{\boldsymbol{\eta}}

\newcommand{\bxi}{\boldsymbol{\xi}}

\newtheorem{theo}{Theorem}[section]
\newtheorem{lmm}[theo]{Lemma}
\newtheorem{df}[theo]{Definition}
\newtheorem{prop}[theo]{Proposition}
\newtheorem{cor}[theo]{Corollary}
\newtheorem{rem}[theo]{Remark}
\newcommand{\Proof}{\noindent {\it Proof: }}
\def\qed{\hfill$\Box$\par\medskip\par\relax}

\begin{document}

\title{A reverse duality for the ASEP with open boundaries}
\author{ 
\and G.M.~Sch\"utz 
}

\maketitle

{\small
Departamento de Matem\'atica, Instituto Superior T\'ecnico,
Universidade de Lisboa,
Av. Rovisco Pais,
1049-001 Lisbon,
Portugal


}

\begin{abstract}
We prove a duality between the asymmetric simple exclusion process 
(ASEP) with non-conservative open boundary conditions and an asymmetric 
exclusion process with particle-dependent hopping rates and conservative 
reflecting boundaries. This is a reverse duality in the sense that the duality 
function relates the measures of the dual processes rather than expectations. Specifically, for a certain parameter manifold of the boundary parameters of the open ASEP this duality expresses the time evolution of a family of shock 
product measures with $N$ microscopic shocks in terms of the time evolution of $N$ particles in the dual process. The reverse duality also elucidates 
some so far poorly understood properties of the stationary matrix 
product measures of the open ASEP given by finite-dimensional matrices. 
\\[.3cm]\textbf{AMS 2010 subject classifications:} 82C22. Secondary: 82C23, 82C26, 60K35
\end{abstract}

\newpage

\section{Introduction}
\setcounter{equation}{0}

Duality expresses time-dependent expectations for one Markov process in 
terms of time-dependent expectations of another Markov process and thus 
is an important tool for the study of stochastic interacting particle systems 
\cite{Ligg85,Ligg99,Jans14}, particularly
if the dual process has some simple properties that are ``hidden'' in
the original process. The key object is the duality function $D$: 
Processes $\bfeta(t)$ with state space $\Omega$ and $\bxi(t)$ with 
state space $\Xi$ respectively are said to be dual w.r.t. a bounded measurable
function $D: \Omega \times \Xi \to \R$
if and only if for all $\bfeta \in \Omega, \bxi \in \Xi$ and $t \geq 0$ one
has
$\E_{\bfeta} D(\bfeta(t),\bxi) = \E_{\bxi} D(\bfeta,\bxi(t))$.
The subscript denotes expectation with the initial states $\bfeta = 
\bfeta(0)$ and $\bxi = \bxi(0)$ resp. 
Considering the duality function 
$D(\bfeta,\bxi)$ to represent a family of functions $f^{\bxi}(\bfeta)
=D(\bfeta,\bxi)$ 
indexed by $\bxi$ then duality expresses for all members of this family the 
expectation of this function at time
$t$ in terms of the expectation at time $t$ of the family of functions
$g^{\bfeta}(\bxi)=D(\bfeta,\bxi)$ indexed by $\bfeta$ of the
dual process.

In the last decade a very considerable body of work on applications of duality 
to interacting particle systems defined on the one-dimensional integer lattice
has emerged see e.g. 
\cite{Beli18,Kuan18,Lin19,Boro20,Chen20,Cari21,Kuan21,Kuan22} 
and references therein specifically 
for very recent work on duality for processes closely related to the ASEP
which is studied in the present paper. To push this field further we introduce 
the concept of reverse duality which expresses the time evolution of 
probability measures of one process in terms of the dual process rather than 
relating time-dependent expectations. Reverse duality can, however,
be rephrased as a conventional duality for the time-reversed processes
and thus extends the framework of conventional duality.
We use it to derive a duality that has so far remained elusive, viz.,
for the asymmetric simple exclusion process (ASEP) on a
one-dimensional integer lattice of $L$ sites with two non-conservative open 
boundaries \cite{Ligg76,Krug91,Derr93a,Schu93b}. On a certain manifold of 
boundary rates for injection and absorption the open ASEP is shown to be 
reverse-dual to an ASEP with particle-dependent jump rates on the same 
lattice but with two conservative reflecting boundaries. Thus we go beyond a
recent proof of self-duality for one open and one reflective boundary condition 
\cite{Kuan21,Kuan22} and an earlier failure to construct a non-trivial duality 
for two open boundaries \cite{Ohku17}.

The reverse duality for the open ASEP constructed here
allows for expressing the time evolution of a certain family of
Bernoulli shock measures with $N$ shocks in terms of the evolution
of the probability measure of the dual ASEP with $N$ particles 
and thus reduces the complexity of the time evolution of the system
with $2^L$ states to $L\choose N$ states. Specifically,
the time evolution of a single shock in the many-particle ASEP
turns into the simple problem of a single biased random walker.
The evolution of initial measures with more than one shock can be 
studied through the complete integrability \cite{Inam94,deVe94} which 
provides an explicit system of Bethe equations for the spectrum of the 
generator \cite{Nepo03b,Nepo04,deGi05,deGi06,Simo09} and for which
results on the eigenfunctions are available \cite{Cram10}. Along 
the way one also obtains insight into some of the properties of the 
stationary matrix product measures of the open ASEP introduced 
in \cite{Derr93a} which on the parameter manifold studied here 
is given by finite-dimensional matrices \cite{Mall97}. 

The paper is organized as follows. In Sec. 2 we introduce reverse
duality in general terms, and define the ASEP with open boundaries, 
shock measures for the ASEP (which give rise to the duality function), 
and the dual shock ASEP with particle-dependent hopping rates 
(which turns out to the reverse dual). The results, in particular the
reverse duality, are stated in precise form as theorems in Sec. 3, with 
comments on the significance of these results. The theorems are 
formally proved in Sec. 4. In Appendix A we show how the generator 
of the open ASEP in the parametrization of the present work is related 
to the integrable Heisenberg quantum spin chain with non-hermitian 
boundary fields of \cite{Nepo03b,Nepo04} and in Appendix B we
discuss some properties of the boundary parameters.

\section{Definitions}
\setcounter{equation}{0}

\subsection{General definitions and conventions}

We consider the finite integer lattice $\Lambda_L = 
[L_-,L_+] \cap \Z$ with $L=L_+-L_-+1 \geq 2$ sites. 
The Kronecker-$\delta$ is defined by
\bel{Def:Kron}
\delta_{x,y} := \left\{ \ba{ll} 1 & \mbox{ if } x=y \\ 
0 & \mbox{ else } \ea \right.
\ee
for $x,y$ from any set. 

For a finite or countably infinite vector space column vectors 
with components $\phi_i$ are denoted by the quantum mechanical 
ket symbol $\ket{\phi}$ while row vectors denoted by the bra symbol
$\bra{\phi}$ have complex conjugate components $\bar{\phi}_i$.
The standard scalar product is written in the Dirac convention
\bel{scp}
\inprod{\phi}{\psi} = \sum_i \bar{\phi}_i \psi_i
\ee
Notice the complex conjugation of the left argument.

\begin{df}[Kronecker product]
\label{Def:Kronecker}
Let $A$ and $B$ be two finite-dimensional matrices with $m_A\geq 1$ 
($m_B\geq 1$) rows and $n_A\geq 1$ ($n_B\geq 1$) columns with matrix 
elements $A_{ij}$ and $B_{kl}$ respectively. The Kronecker product
$A\otimes B$ is a $m_A m_B \times n_An_B$-matrix $C$ with matrix elements
\bel{Kronecker}
C_{(i-1)m_B+k,(j-1)n_B+l}=A_{ij} B_{kl}
\ee 
with $i \in \{1,\dots,m_A\}, j \in \{1,\dots,n_A\}, k \in \{1,\dots,m_B\}, 
l \in \{1,\dots,n_B\}$. The $n$-fold Kronecker product of a matrix $A$
with itself is denoted by $A^{\otimes n} $ and by convention 
$A^{\otimes 0} = 1$ and $A^{\otimes 1} = A$.
\end{df}

\begin{df}[Local and bond operators]
(i) With the two-dimensional unit matrix denoted by 
$\mathds{1}$ and an arbitrary two-dimensional matrix $\mathsf{a}$
the Kronecker product
\be 
\mathsf{a}_k := \mathds{1}^{\otimes (k-L_-)} \otimes \mathsf{a} 
\otimes \mathds{1}^{\otimes (L_+-k)}, \quad k \in \Lambda_L
\ee
is called a local operator acting on site $k\in\Lambda_L$.\\
\noindent (ii) For any four-dimensional matrix $\mathsf{b}$, bond 
operators $\mathsf{b}_{k,k+1}$ with $L_- \leq k \leq L_+-1$ are 
defined by the Kronecker product
\be 
\mathsf{b}_{k,k+1} = \mathds{1}^{\otimes (k-L_-)} \otimes 
\mathsf{b} \otimes \mathds{1}^{\otimes (L_+-1-k)}.
\ee
\end{df}

\subsection{Duality and reverse duality}

For many applications it is convenient to express the generator of a 
Markov process $\bfeta(t)$ with countable state space in terms of 
the intensity matrix $W$ whose off-diagonal matrix elements 
$W_{\bfeta,\bfeta'}$ are the transition rates from configuration $\bfeta$ 
to $\bfeta'$ and with diagonal elements given by the negative sum 
$W_{\bfeta\bfeta} = - \sum_{\bfeta'} W_{\bfeta,\bfeta'}$ of outgoing 
transition rates from state $\bfeta$. In terms of intensity matrices $
W$ for the process $\bfeta(t)$ and $Q$ for the process $\bxi(t)$, 
and with the duality matrix $D$ with the matrix elements 
$D_{\bfeta\xi} = D(\bfeta,\bxi)$ duality is expressed in matrix form 
as the relation \cite{Sudb95}
\bel{duality}
W D = D Q^T. 
\ee
This notion of duality has the following natural generalization.

\begin{df}
\label{Def:revduality}
Let $\bfeta(t)$ be a Markov process with intensity matrix $W$. A 
Markov process $\bxi(t)$ with intensity matrix $Q$ satisfying 
\bel{revdualitydef}
R W = Q^T R
\ee
with a duality matrix $R$ with matrix elements $R_{\bxi\bfeta}$ 
is called reverse dual to $\bfeta(t)$ w.r.t. the duality function 
$R(\bxi,\bfeta)=R_{\bxi\bfeta}$.
\end{df}

Obviously, reverse duality is a mutual relation between two processes
as can be seen by transposing both sides of the defining equation
\eref{revdualitydef} which then says that $\bfeta(t)$ is a reverse
dual to $\bxi(t)$ with duality matrix $R^T$. To give some motivation 
for the notion of revsere duality assume that the processes $\bfeta(t)$ and 
$\bxi(t)$ have invertible invariant measures $\pi_W$ and $\pi_Q$ resp.,
i.e., $\pi_W \neq 0 $ for all $\bfeta$ and $\pi_Q \neq 0$ for all $\bxi$.
This implies that that the diagonal matrices $\hat{\pi}_W^{-1}$ with 
diagonal entries $\pi_W^{-1}(\bfeta)$ and $\hat{\pi}_Q^{-1}$ with 
diagonal entries $\pi_Q^{-1}(\bxi)$ exist. By definition, the intensity 
matrices for the reversed processes are then given by $W_{rev} = 
\hat{\pi}_W^{-1} W^T \hat{\pi}_W$ and $Q_{rev} = 
\hat{\pi}_Q^{-1} Q^T \hat{\pi}_Q$, resp., so that $W^T = 
\hat{\pi}_W W_{rev}  \hat{\pi}_W^{-1}$ and 
$Q = \hat{\pi}_Q^{-1} Q_{rev}^T  \hat{\pi}_Q$. 
After transposition and with the duality matrix $D=\hat{\pi}_W^{-1} R^T  
\hat{\pi}_Q^{-1}$ the defining relation \eref{revdualitydef} for reverse 
duality turns into the conventional duality relation $W_{rev} D =  D 
Q_{rev}^T$ for the reversed processes. In other words, for reversible 
processes reverse duality implies a conventional duality for the 
reversed processes and vice versa.

While conventional duality relates expectation values of the dual processes, 
reverse duality can be used to relate the measures
for all times $t\geq 0$. To see this, we denote by $\mu(t)$ 
the probability of the configuration $\bfeta(t)$ for some unspecified initial measure
$\mu(0)$ and for the process $\bxi(t)$ we denote by $P(\bxi,t|\bxi',0)$ the
transition probability at time $t$ from a configuration $\bxi'$ to a configuration 
$\bxi$.

\begin{theo}
\label{Theo:measure1}
For a countable set $\Xi$ let $\mathcal{F}_{\Xi}$ be a family of probability 
measures for the process $\bfeta(t)$ with countable state space $\Omega$
and let the process $\bxi(t)$ with state space $\Xi$ be reverse dual to 
$\bfeta(t)$ with duality function $R(\bxi,\bfeta) = \mu^{\bxi}_{\bfeta} 
\in \mathcal{F}_{\Xi}$. Then at time $t\geq 0$ the measure $\mu^{\bxi}(t)$
is given by the convex combination
\bel{measure1} 
\mu^{\bxi}(t) =  
\sum_{\bxi'} P(\bxi,t|\bxi',0)  \mu^{\bxi'}(0)
\ee
of initial measures $\mu^{\bxi'}_0\in\mathcal{F}_{\Xi}$.
\end{theo}

\Proof The time evolution property \eref{measure1} of the measure follows 
directly from the semigroup property of the Markov processes $\bfeta(t)$ 
and $\bxi(t)$ and Definition \ref{Def:revduality} of reverse duality. 
\qed

Thus if the duality function $R(\bxi,\bfeta)$ represents a family
of probability measures $\mu^{\bxi}_{\bfeta}$ indexed by $\bxi$ 
then reverse duality expresses the Markovian time evolution of all 
members of this family of measures in terms of the probability 
measure of the reverse dual process starting from configuation 
$\bxi$.\footnote{Duality 
between measures is discussed in a different setting also in 
\cite{Chun16}.} In contrast, considering the duality function of conventional duality to be member of a measurable family of functions $f^{\bxi}(\bfeta)
:=D(\bfeta,\bxi)$ of the states $\bfeta$ and denoting the expecation of this
function at time $t\geq 0$ w.r.t. some unspecified initial measure 
by $\E f^{\bxi}(t)$, conventional duality reads
\be
\E f^{\bxi}(t) = \sum_{\bxi'}  P(\bxi,t|\bxi',0) \E f^{\bxi'}(0). 
\ee
This highlights the difference between reverse duality and conventional
duality in terms of the role of the transition probability of the (reverse) 
dual process.

For reversible processes
$\bxi(t)$ with reversible measure $\pi$ Theorem \ref{Theo:measure1} 
can be stated in a slightly modified form.

\begin{theo}
\label{Theo:measure2}
For a countable set $\Xi$ let $\mathcal{F}_{\Xi}$ be a family of probability measures for the process $\bfeta(t)$ with countable state space $\Omega$
and let the process $\bxi(t)$ with state space $\Xi$ be (a) reversible with reversible measure $\pi$ and (b) reverse dual to $\bfeta(t)$ with duality function
$R(\bxi,\bfeta) = \pi(\bxi) \mu^{\bxi}_{\bfeta} \in \mathcal{F}_{\Xi}$. Then at time $t\geq 0$ the measure $\mu^{\bxi}(t)$
is given by the convex combination
\bel{measure2}
\mu^{\bxi}(t) =  
\sum_{\bxi'} P(\bxi',t|\bxi,0)  \mu^{\bxi'}(0)
\ee
of initial measures $\mu^{\bxi'}(0)\in\mathcal{F}_{\Xi}$.
\end{theo}

\Proof With 
$S := \hat{\pi}^{-1} R$
reverse duality as defined in \eref{revdualitydef} implies for $Q=Q^{rev}$ the intertwining relation
\bel{intdualitydef}
S W = Q  S.
\ee
The time evolution
property \eref{measure2} then 
follows from the semigroup property of the Markov
processes $\bfeta(t)$ and $\bxi(t)$ and reversibility of $\bxi(t)$.
\qed

\begin{rem}
The difference to Theorem \ref{Theo:measure1} appears in the
time-reversal of the transition probability of the reverse dual process
and in the relationship between the duality function $R(\bxi,\bfeta)$ and the 
measures $\mu^{\bxi}$.
\end{rem}

\begin{rem}
For self-dual processes with $W=Q$ the intertwining relation 
\eref{intdualitydef} expresses a symmetry of the generator 
\cite{Schu94,Giar09,Cari19} and for invertible
$S$ it constitutes a similarity transformation between processes
with intensity matrices $W,Q$ \cite{Henk95}.
\end{rem}

\subsection{Intensity matrices and quantum Hamiltonians}

We recall that by construction the sum of all matrix elements in each row 
of the intensity matrices add up to zero. It is well-known that the 
transposed intensity matrix often takes the form of a Hamiltonian 
operator $H= - W^T$ for quantum mechanical spin chains 
\cite{Alca94,Lloy96,Schu01} and concepts from quantum mechanics can be 
imported to solve probabilistic problems even though in general $H$ is not 
a Hermitian matrix. To illustrate the point in the context of duality we 
mention that an invariant measure is a ground state vector of the quantum
system and non-trivial self-dualities arise from non-Abelian symmetries 
of the generator \cite{Schu95,Schu97,Giar09,Cari19}. An general relation 
between duality and spectral properties of the generator - which play a crucial
role in quantum mechanics and can be computed exactly for integrable systems -
was recently worked out in \cite{Redi18}. 
Complete integrability can also be used to directly
construct novel dualities \cite{Schu20}. 

With this in mind we express below duality and reverse duality
in terms of the negative transpose $H$ of the intensity matrix
and in slight abuse of language we also call $H$ the
generator of the process. Time-invariance of a measure $\pi$ 
means that the vector $\ket{\pi}$ is a right eigenvector of $H$ 
with eigenvalue $0$.
Defining the diagonal matrix $\hat{\pi}$ with the components $\pi(\cdot)$
on the diagonal, the generator of the time-reversed process is given by
\bel{revdef}
H_{rev} := \hat{\pi} H^T \hat{\pi}^{-1} 
\ee
and has the same invariant measure $\pi$.
With $H=-W^T$ and $G=-Q^T$ duality as defined in \eref{duality}
then reads
\bel{dualityH}
D^T H = G^T D^T 
\ee
with transposed duality matrix $D^T$ with matrix elements
$D^T_{\bxi\bfeta} = D(\bfeta,\bxi)$.
Reversed duality as defined in \eref{revdualitydef} becomes
\bel{revdualityH}
H R^T = R^T G^T
\ee
and for reversible $G$ this yields the intertwining relation
\bel{intdualityH}
H S^T = S^T G 
\ee
with $S^T=R^T \hat{\pi}^{-1}$.

To rephrase Theorem \ref{Theo:measure2} we use quantum mechanical bra-ket notation for vectors and the scalar product \eref{scp}. Then a probability measure for the process
$\bfeta(t)$ is a ket vector
$\ket{\mu(t)}$ whose components $\mu_{\bfeta}(t)$ are the 
probabilities of configuration $\bfeta(t)$. The Markov semigroup property
yields $\ket{\mu(t)} = [\bra{\mu(0)} \exp{(Wt)}]^T
= \exp{(-Ht)} \ket{\mu(0)}$ so that the 
transition probability $\mu(\bfeta',t|\bfeta,0)$ from a configuration $\bfeta$ to a configuration 
$\bfeta'$ is given by
\be 
\mu(\bfeta',t|\bfeta,0) = \bra{\bfeta} \rme^{Wt} \ket{\bfeta'} =
\bra{\bfeta'} \rme^{-Ht} \ket{\bfeta}.
\ee
and similarly for the process $\bxi(t)$ generated by $G$. The evolution
formula \eref{measure2} of Theorem \ref{Theo:measure2} reads in 
vector form
\bel{measure3}
\ket{\mu^{\bxi}(t)} =  
\sum_{\bxi'} P(\bxi',t|\bxi,0)  \ket{\mu^{\bxi'}(0)}.
\ee
with the transition probability $P(\bxi',t|\bxi,0)$ of the dual process.

\begin{rem}
Under the conditions of Theorem \ref{Theo:measure2} the family
$\mathcal{F}_\Xi$ of probability vectors belongs to a right-invariant subspace of the generator $H$ with dimension $d_\Xi \leq |\mathcal{F}_\Xi|$.
\end{rem}

\subsection{ASEP with open boundaries}

\subsubsection{Generator}

In the ASEP each lattice site $k\in\Lambda_L$ can be occupied or empty as 
expressed by the occupation numbers $\eta_k\in\{0,1\}$. The particle 
configuration on the whole lattice $\lambda_L$ is denoted by $\bfeta := 
\{\eta_{L_-}, \eta_{L_-+1}, \dots, \eta_{L_+}\}$. For a given 
configuration $\bfeta$ we also define
the flipped configuration $\bfeta^{k}$ with occupation numbers
\bel{etaflipdef}
\eta^{k}_l = \left\{\ba{ll}
1-\eta_{k} & \mbox{if } l = k \\
\eta_{l} & \mbox{else} .
\ea \right.
\ee
and the
swapped configuration $\bfeta^{kk+1}$ with occupation numbers
\bel{etaswapdef}
\eta^{kk+1}_l = \left\{\ba{ll}
\eta_{k+1} & \mbox{if } l = k \\
\eta_{k} & \mbox{if } l = k+1 \\
\eta_{l} & \mbox{else} .
\ea \right.
\ee
Expectations are denoted by $\exval{\cdot}$. In particular, 
$\exval{\eta_k}$ is called the local density.

Informally, the Markovian dynamics of the open ASEP can be described as follows. A particle on a bulk site $k \in [L_-+1,L_+-1]$ jumps
with rate $r$ ($\ell$) independently of the other particles
to right (left) neighbouring site $k+1$ ($k-1$), provided
the target site is empty. Otherwise the jump attempt is rejected.
On the left boundary site $L_-$ a particle jumps with rate $r$ to the
right (if site $L_-+1$ is empty) or is annihilated with rate $\gamma$.
If site $L_-$ is empty then a particle is created with rate $\alpha$.
On the right boundary site $L_+$ a particle jumps with rate $\ell$ to the
left (if site $L_+-1$ is empty) or is annihilated with rate $\beta$.
If site $L_+$ is empty then a particle is created with rate $\delta$.

With the local jump rate
\bel{ASEPjumprate}
w_k(\bfeta) = r \eta_{k} (1-\eta_{k+1}) + \ell (1-\eta_{k}) \eta_{k+1}
\ee
and the boundary rates
\bel{ASEPratebl}
b^-(\bfeta) = \alpha (1-\eta_{L_-}) + \gamma \eta_{L_-}, \quad
b^+(\bfeta) = \delta (1-\eta_{L_+}) + \beta \eta_{L_+}
\ee
the generator $\mathcal{L}$ of this process is thus given by
\bel{genASEP}
\mathcal{L} = \mathcal{L}^{bulk} + \mathcal{L}^{-} + \mathcal{L}^{+}
\ee
with the bulk generator
\bel{genASEPbulk}
\mathcal{L}^{bulk} f(\bfeta) = \sum_{k=L_-}^{L_+-1}
w_k(\bfeta) [f(\bfeta^{kk+1}) -  f(\bfeta)]
\ee
and the boundary generators
\bel{genASEPboundary}
\mathcal{L}^{-} f(\bfeta) = b^-(\bfeta) [f(\bfeta^{L_-}) -  f(\bfeta)], \quad
\mathcal{L}^{+} f(\bfeta) = b^+(\bfeta) [f(\bfeta^{L_+}) -  f(\bfeta)].
\ee
Throughout this work it is assumed that all rates $r,\ell,\alpha,\beta,\gamma,\delta$ 
are strictly positive and that $r \neq \ell$.

To write the corresponding intensity matrix we denote the
two-dimensional unit matrix by $\mathds{1}$ and define the matrices
\bea
w^b & = & \left(\ba{cccc}
0 & 0 & 0 & 0 \\
0 & - \ell & \ell & 0 \\
0 & r & - r & 0 \\
0 & 0 & 0 & 0 
\ea \right) \\
w^{-} & = & \left(\ba{cc}
- \alpha & \alpha \\
\gamma & - \gamma
\ea \right) \\
w^{+} & = & \left(\ba{cc}
- \delta & \delta \\
\beta & - \beta
\ea \right) .
\eea
and the Kronecker products $w^{-}_{L_-} := w^{-} \otimes \mathds{1}^{\otimes L_+ - L_-}$, $w^{+}_{L_+} := \mathds{1}^{\otimes L_+ - L_-} \otimes w^{+}$, and
$w^b_{k,k+1} := \mathds{1}^{\otimes -L_-+k} \otimes w^b \otimes \mathds{1}^{\otimes L_+ -1 - k}$. In terms of these matrices
the intensity matrix for the open ASEP is given by
\bel{MASEP}
W = W_{bulk} + W_{edge}
\ee
with
\be 
W_{bulk} = \sum_{k=L_-}^{L_+-1} w^b_{k,k+1}, \quad
W_{edge} = w^{-}_{L_-} + w^{+}_{L_+}.
\ee
In $H:= - W^T$ one recognizes a non-hermitian version of the
XXZ quantum spin chain with non-diagonal boundary fields \cite{Nepo03b,Nepo04},
see the appendix.

To make contact with the well-known properties of the XXZ chain, the computations below are done using
$H$ rather than $W$. Defining the particle number operator
\be 
\hat{n} := \left(\ba{cc}
0 & 0 \\
0 & 1
\ea \right) 
\ee 
and the Kronecker product 
$\hat{n}_{k} := \mathds{1}^{\otimes -L_-+k} \otimes \hat{n} \otimes 
\mathds{1}^{\otimes L_+ - k}$ 
one gets
\bel{generator}
H = H_{bulk} + H_{edge}
\ee
with
\bel{Hbulk}
H_{bulk} = \sum_{k=1}^{L-1} \tilde{h}_{k,k+1}, \quad 
H_{edge} = \tilde{h}^{-}_{L_-} + \tilde{h}^{+}_{L_+}
\ee
where
\bea
\tilde{h}_{k,k+1} & = & - (w^b_{k,k+1})^T + (r-\ell) \left(\hat{n}_{k+1} - \hat{n}_{k}\right)
\, = \, \left(\ba{cccc}
0 & 0 & 0 & 0 \\
0 & r & - r & 0 \\
0 & - \ell & \ell & 0 \\
0 & 0 & 0 & 0 
\ea \right)_{k,k+1} ,
\label{hbulk} \\
\tilde{h}^{-}_{L_-} & = & - (w^-_{L_-})^T   + (r-\ell) \hat{n}_{L_-}
\, = \, \left(\ba{cc}
\alpha & - \gamma \\
- \alpha & \gamma + r - \ell
\ea \right)_{L_-} ,
\label{hbm} \\
\tilde{h}^{+}_{L_+} & = & - (w^+_{L_+})^T + (r-\ell) \hat{n}_{L_+}
\, = \, \left(\ba{cc}
\delta & - \beta \\
- \delta & \beta - (r - \ell)
\ea \right)_{L_+} .
\label{hbp}
\eea
Notice the addition of the discrete gradient $(r-\ell) \left(\hat{n}_{k+1} - \hat{n}_{k}\right)$ to the local bulk generators which cancel in \eref{generator} with the
corresponding contributions to the boundary generators in \eref{hbm}
and \eref{hbp}.

\subsubsection{Particle current}

Particle conservation in the bulk is reflected in the discrete continuity
equation
\be 
\mathcal{L} \eta_k = j_{k-1} - j_k
\ee 
where
\bel{instcurrbulk}
j_k = r \eta_k(1-\eta_{k+1}) - \ell (1-\eta_k) \eta_{k+1}, \quad
L_- \leq k < L_+
\ee
is the instantaneous current across bond $(k,k+1)$.
At the boundaries one has
\be 
\mathcal{L} \eta_{L_-} = j^- - j_{L_-}, \quad
\mathcal{L} \eta_{L_+} = j_{L_+-1} - j^+
\ee 
where
\be 
j^- = \alpha(1-\eta_{L_-}) - \gamma \eta_{L_-}, \quad
j^+ = \beta \eta_{L_+} - \delta(1-\eta_{L_-}) 
\ee
are the boundary currents with a source term a sink term each.
Stationarity implies constant expectation
\bel{jstat}
j := \exval{j^-} = \exval{j_k} = \exval{j^+}
\ee
in the invariant measure which was first computed for the TASEP ($\gamma=\delta=\ell=0$) in two different ways, viz., with a matrix product ansatz in \cite{Derr93a} and simultaneously by generating function techniques in \cite{Schu93b}. The general case was solved
subsequently in \cite{Sand94b}.

\subsubsection{Invariant matrix product measures}

To discuss some properties of the invariant measure it is convenient to
introduce the jump asymmetry $q$ and the time scale $w$ by
\be
q := \sqrt{\frac{r}{\ell}}, \quad w:= \sqrt{r\ell}
\label{qdef} 
\ee
and to parametrize the boundary rates as
\bea
& & \alpha = (r+\omega_-) \rho_-, \quad \gamma = (\ell+\omega_-) (1-\rho_-)
\label{lbr} \\
& & \beta = (r+\omega_+) (1-\rho_+), \quad \delta = (\ell+\omega_+)\rho_+
\label{rbr}
\eea
where the parameters $\omega_\pm$ may be interpreted as
boundary jump barriers and, as seen below, the parameters 
$\rho_\pm$ play the role of boundary densities.
It is convenient to also consider the fugacity
\bel{fugadef}
z \equiv  z(\rho) = \frac{\rho}{1-\rho}
\ee
as a function of the density $\rho$. For an indexed density $\rho_i$ 
with $i$ from any index set we shall use the notation
\be 
z_i \equiv z(\rho_i)
\ee
as shorthand. 

The invariant measure of the open ASEP is conveniently characterized in terms of the functions \cite{Sand94b}
\be
\kappa_\pm(x,y) := \frac{1}{2x} (y - x + r - \ell \pm
\sqrt{(y - x + r - \ell))^2 + 4x y}) 
\label{kappadef}
\ee
which are the roots of the quadratic equation
\be  
x \kappa^2 - (y - x + r - \ell) \kappa - y = 0 
\ee 
and which are related by
\be 
\kappa_+(x,y)\kappa_-(x,y) = - \frac{y}{x}.
\label{kpm}
\ee

As shown in \cite{Mall97} on the parameter manifold $\mathcal{B}_N$
defined by the relation
\be
\label{bNshockcond} 
\kappa_+(\alpha,\gamma) \kappa_+(\beta,\delta) = q^{2N} 
\ee
the invariant measure for the open ASEP
can be expressed in terms of a matrix 
product measure (MPM) with $(N+1)$-dimensional representation matrices of a quadratic algebra satisfied by the
matrices of the matrix product ansatz developed in \cite{Derr93a}.
On the submanifold $\mathcal{B}^M_N\subset\mathcal{B}_N$
defined by the
further constraint
\bel{bNshockcond2}
\kappa_-(\alpha,\gamma) \kappa_-(\beta,\delta) = q^{-2M}, \quad
1 \leq M \leq N
\ee
the MPM has the special property that it exists only for system
sizes $L > N - M +1$ \cite{Mall97} whereas away from this submanifold
it exists for all system sizes \cite{Essl96}. 
Boundary parameters satisfying \eref{bNshockcond} and \eref{bNshockcond2} 
for $1\leq N \leq L$ define the setting that the main results
proved below apply to.

To get further insight into the significance of the parameter manifold
$\mathcal{B}_N$ we note that for $N=0$ the condition \eref{bNshockcond} 
reduces to $z_+=z_-$. As shown in \cite{Sand94b}, with this constraint 
on the boundary parameters the invariant measure is a Bernoulli product 
measure with density $\rho_+=\rho_-$ and the manifolds \eref{aga} and 
\eref{bdb} in Appendix \ref{App:B} turn into the stationarity
condition \eref{jstat} for the expectation of the current. Following 
\cite{Ligg76} this observation may be understood as follows. In a semi-infinite 
system with $L_+=+\infty$ the left boundary rates \eref{lbr} allow for an 
invariant Bernoulli product measure with stationary density $\rho_- \leq 1/2$ 
and stationary current given by \eref{aga} (corresponding to $j=\exval{j^-}$) 
while in a semi-infinite system with $L_-=- \infty$ the right boundary
rates \eref{rbr} allow for an invariant Bernoulli product measure
with a stationary density $\rho_+ \geq 1/2$ and stationary
current given by \eref{bdb} (corresponding to $j=\exval{j^+}$). 
When $\rho_+=\rho_-$ or $\rho_+=1-\rho_-$ the two currents $\exval{j^-}$ 
and $\exval{j^+}$ match. Then stationarity of the Bernoulli product
measure holds even in the finite lattice $\Lambda_L$
for any $L_\pm$ and any density $\rho:=\rho_+=\rho_-$.

On the other hand, when $\rho_+=1-\rho_-$ stationary can only be attained
with a non-homogeneous density profile which has been argued to be a shock 
profile \cite{Schu93b,Schu97b} analogous to a domain wall in driven 
non-equilibrium systems \cite{Evan98b,Lahi00,Clin03,Kafr03,Chak16} but with 
diffusive fluctuations of the domain wall position \cite{Dudz00,Sant02,deGi06}
proved rigorously for the ASEP on the infinite integer lattice in \cite{Ferr94}.
Any mismatch of the boundary currents in the range
$\rho_- < 1/2$ and $\rho_+ > 1/2$ is then expected to be still realized
by a shock profile, but with a shock localized close to the left boundary
for $\exval{j^+}<\exval{j^-}$ or close to the right boundary for
$\exval{j^+}>\exval{j^-}$, thus explaining the
origin of the phase diagram of the open ASEP in terms of domain wall
dynamics \cite{Kolo98b}. From a mathematical perspective, this
reasoning is a conjecture based on microscopic shock stability
proved only in the hydrodynamic limit \cite{Baha12}. As will be seen below, 
the reverse duality constructed below provides a rigorous proof of this 
conjecture on microscopic scale for the
parameter manifold $\mathcal{B}^1_N$ defined by \eref{bNshockcond} 
and \eref{bNshockcond2}.

Finally we note that for system size $L=N$ the invariant measure on the
parameter manifold $\mathcal{B}^1_N$ is a zero-current
blocking measure restricted to the finite lattice \cite{Bryc19}.
This blocking measure, originally defined for the infinite system \cite{Ligg85}, 
is a product measure with strictly increasing marginal fugacities 
$z_k \propto q^{2(k-L_-)}$, $k\in\Lambda_L$,
is also an invariant measure of the ASEP with reflecting boundaries \cite{Sand94a} where $\alpha=\beta=\gamma=\delta=0$.

\subsection{Shocks in the ASEP}

\subsubsection{Bernoulli shock measures}

In the spirit of \cite{Beli02,Beli18} we extend the notion of Bernoulli 
shock measures to the finite lattice. To this end and to allow for more 
compact notation auxiliary boundary sites are introduced.

\begin{df}
\label{Def:shockmeas} (Bernoulli shock measures)
For auxiliary boundary sites $x_0:=L--1$ and $x_{N+1} := L_++1$ and a nonempty ordered set $\bfx := \{x_1,\dots,x_N)\}$  of $N$ lattice sites 
the product measure 
\bel{shockmeas}
\mu^{\bfx}_{\bfeta} = \prod_{k=L_-}^{L_+} p^{\bfx}_{\eta_k}
\ee
with marginals 
\bel{shockmarg} 
p^{\bfx}_{\eta_k} = \left\{ \ba{ll}
(1-\rho^\star_i) (1-\eta_{k}) + \rho^\star_i \eta_{k}  & k = x_i, \quad 1 \leq i \leq N \\
(1-\rho_i) (1-\eta_k) + \rho_i \eta_k  & x_i < k < x_{i+1}, \quad 0 \leq i \leq N \\
\ea\right. 
\ee
is called a Bernoulli shock measure with $N$ microscopic shocks at positions $x_i \in \Lambda_L$ and left boundary density $\rho_0$, right boundary density $\rho_N$,  bulk densities $\rho_i$ for $1\leq i \leq N-1$, and shock densities $\rho^\star_i$ for $1\leq i \leq N$. 
\end{df}
This definition is motivated by the fact that seen from a shock position
the local particle configuration of an infinite system becomes uncorrelated 
at large distances \cite{Ferr91,Derr93}. The Bernoulli shock measures
extend this property to all lattice sites which was found in \cite{Derr97}
for specific shock densities whose fugacities satisfy the relation 
\be
\label{shockcondi}
\frac{z_i}{z_{i-1}} = q^{2} 
\ee 
and second-class particles at the shock positions.

The vector representation of a shock measure with $N$ shocks
is given by
\be 
\ket{\mu^\bfx} :=  \sum_{\bfeta} \mu^{\bfx}_{\bfeta} \ket{\bfeta}
\ee
and with the column vector $\ket{\rho_{\cdot}} := (1-\rho_{\cdot},\rho_{\cdot})^T$ 
one gets the Kronecker product
\bel{shockmeasvec}
\ket{\mu^\bfx} = \ket{\rho_{0}}^{\otimes (x_1 - L_-)} \otimes \ket{\rho^{\star}_1} \otimes \ket{\rho_{1}}^{\otimes (x_2 - x_1 - 1)} \otimes \dots \otimes \ket{\rho^{\star}_N} \otimes \ket{\rho_{N}}^{\otimes (L_+ - x_N)}
\ee 
which reflects the factorized form of the shock measure.

\subsubsection{Microscopic and macoscopic shocks}

For a Bernoulli shock measure the expectation of the particle current 
\eref{instcurrbulk} inside a domain with density $\rho_i$ is given by
\be
j_i \equiv  j(\rho_i) = (r-\ell) \rho_i(1-\rho_i) .
\label{ASEPcurrent} 
\ee
As shown in \cite{Beli02} for shock densities 
$\rho^{\star}_i = 1$ the quantities
\bea 
d^{\ell}_{i} & = & 
\frac{j_{i-1}}{\rho_{i}-\rho_{i-1}} 
\label{dildef} \\
d^{r}_{i} & = & 
\frac{j_{i}}{\rho_{i}-\rho_{i-1}} 
\label{dirdef} 
\eea
are the shock hopping rates of a microscopically stable shock with 
bulk densities $\rho_{i-1},\rho_{i}$ satisfying \eref{shockcondi} which
we therefore call the
microscopic shock stability condition.
In \cite{Kreb03,Bala10,Beli15b,Beli18} the same microscopic shock jump 
rates were found with a second-class particle \cite{Ferr91} at the shock 
position which is equivalent to $\rho^{\star}_i = 0$. 
The average shock velocity of  shock $i$ 
between bulk densities $\rho_i$ and $\rho_{i+1}$
is the difference 
\bel{shockvel}
v_i = d^{r}_{i} - d^{\ell}_{i}
\ee 
of shock jump rates
and the variance of the shock position is described by the diffusion
coefficient
\bel{shockdiff}
D_i = \half (d^{r}_{i} + d^{\ell}_{i})
\ee 
which appear also on macroscopic scale in the random motion of a 
macroscopically stable shock in the ASEP \cite{Ferr91,Ferr94}. Notice, 
however, that because of the convexity of the current-density relation 
$j\propto\rho(1-\rho)$ only the milder condition Rankine Hugoniot condition 
\cite{Reza95} is required for the macroscopic stability of a shock 
which for 
the ASEP only implies $\rho_{i+1} > \rho_i$ rather than the strict condition 
\eref{shockcondi}. For consecutive shocks this convexity also
yields $v_i > v_{i+1}$ for all $i$
so that after finite microscopic time consecutive shocks become very close and form a bound state for all subsequent times \cite{Beli02}. On macroscopic scale this phenomenon corresponds to a 
coalescence of shocks \cite{Ferr00}. Also other particle-conserving 
models with random-walking shocks are known \cite{Kreb03,Bala04,Jafa09,Beli13,Bala19}

\subsubsection{Shock random walks and MPMs for the open ASEP}

A shock measure with $N$ stable shocks which all satisfy the 
microscopic stability criterion \eref{shockcondi} 
then has the property
\be
\label{Nshockcond}
\frac{z_N}{z_{0}} = q^{2N}.
\ee 
If the boundary parameters of the open ASEP satisfy \eref{bNshockcond} 
with $z_-=z_0$ then $z_+ = z_N$ then also \eref{Nshockcond} is satisfied 
and the $(N+1)$-dimensional representations of \cite{Mall97} of the
stationary matrix product algebra of \cite{Derr93a} are convex
combinations of homogenous Bernoulli product measures with fugacities $z_0$ and $z_N$ and Bernoulli shock measures with $N$ shocks
satisfying \eref{shockcondi}
and \eref{bNshockcond} \cite{Jafa07}. Therefore we refer to 
\eref{bNshockcond} as boundary shock stability condition.

In fact, it was suggested even earlier in \cite{Kreb03} that the shocks 
whose invariant distribution in the finite system is described by the 
MPM remain stable during the stochastic time evolution
and perform a random walk dynamics. 
Here we focus on boundary rates satisfying also the
condition \eref{Nshockcond} and use reverse duality to prove this conjecture and make it precise. In particular, we identify the
shock densities $\rho^\star_i$ which allows us to elucidate special
properties of MPMs satisfying \eref{bNshockcond2}.

\subsection{Shock ASEP with particle-dependent hopping rates}

Using duality it is proved below that on microscopic scale shocks
in the open ASEP on $\Lambda$, or more precisely the microscopic shock 
positions perform a random motion that we introduced in more 
general form in \cite{Beli18} under the name {\it shock exclusion process} 
in the context of self-duality for a multi-species ASEP. 
Here we go beyond \cite{Beli18} by considering non-conservative
open boundaries instead of the more straightforward case of 
conservative dynamics on $\Z$. However, we focus on only
one species of particles.

For a single species of exclusion particles the shock exclusion process, has 
particles, labelled by $i$ with $1\leq i \leq N$ that are located on lattice
sites $x_i$ with the single-file condition $x_0 < x_1 < x_2 \dots < x_N < 
x_{N+1}$ with the fixed auxiliary boundary coordinates $x_0 := L_--1$
and $x_{N+1} = L_++1$. Thus the configuration of all particles is 
represented by the coordinate vector $\bfx:=(x_1,x_2,\dots,x_{N})$
and each particle has its individual jump rate
\be
w_{i}(\bfx) = w^{\ell}_{i}(\bfx) + w^{r}_{i}(\bfx)
\label{wishockasep}
\ee
with left- and right jump rates
\bea
w^{\ell}_{i}(\bfx) & = & d^{\ell}_i (1-\delta_{x_{i},x_{i-1}+1})
(1-\delta_{x_{i},x_{i+1}}) \label{wil} \\
w^{r}_{i}(\bfx) & = & d^{r}_i (1-\delta_{x_{i},x_{i+1}-1})
(1-\delta_{x_{i},x_{i-1}}) \label{wir}
\eea
given by the shock jump rates \eref{dildef} and \eref{dirdef}.
The terms with the Kroneckersymbols entail the single-file exclusion principle
and the boundary terms $x_0$ and $x_{N+1}$ appearing in 
the rates $w_1(\bfx)$ and $w_N(\bfx)$ express the particle number
conservation at the reflecting boundaries.

The generator of the shock ASEP can be defined more formally by 
introducing for a given coordinate vector $\bfx = (x_1,\dots,x_N)$ 
the locally shifted coordinate vectors
$\bfx^{i\pm}$ with the coordinates 
\be
x^{i\pm}_j = x_j \pm \delta_{i,j}.
\ee

\begin{df} 
\label{Def:shockprocess} (Shock exclusion process) 
The $N$-particle shock exclusion process of particles $i$,  $1\leq i \leq N$ that are located on lattice
sites $x_i \in \Lambda_L$ with the single-file condition $x_0 < x_1 < x_2 \dots < x_N < x_{N+1}$ and with jump rates \eref{wil}, \eref{wir}
 is defined by the generator
\bel{shockgen}
\mathcal{M} f(\bfx) = \sum_{i=1}^N \mathcal{M}_i f(\bfx)
\ee
with the single-particle hopping generators
\be 
\label{shockgenlocal}
\mathcal{M}_i f(\bfx) = 
w_i^{\ell}(\bfx)(f(\bfx^{i-})-f(\bfx)) + w_i^{r}(\bfx)(f(\bfx^{i+})-f(\bfx))
\ee
for fixed parameters
$\rho_i \in [0,1]$, $0\leq i \leq N$.
\end{df}

\begin{rem}
\label{Rem:RW}
For $N=1$ the shock exclusion process reduces to a biased random walk
on $\Lambda_L$ with generator 
\be 
\label{RWgen}
\mathcal{M} f(x) = 
w^{\ell}(x)(f(x-1)-f(x)) + w^{r}(x)(f(x+1)-f(x))
\ee
and
with left- and right jump rates
\bea
w^{\ell}(x) & = & d^{\ell}_1 (1-\delta_{x,L_-})
(1-\delta_{x,L_++1}) \label{wil1} \\
w^{r}(x) & = & d^{r}_1 (1-\delta_{x,L_+})
(1-\delta_{x,L_--1}) \label{wir1} 
\eea
for any fixed pair of densities $\rho_0,\rho_1\in[0,1]$. The intensity
matrix $Q$
is the tridiagonal Toeplitz matrix with matrix elements 
\bea 
Q_{xy} & = & d^{\ell}_1 (1-\delta_{x,L_-}) (\delta_{y,x-1} - \delta_{y,x})
+ d^{r}_1 (1-\delta_{x,L_+}) (\delta_{y,x+1} - \delta_{y,x}) 
\label{intrw}
\eea
for $x,y\in\Lambda_L$.
\end{rem}

\section{Reverse duality for the open ASEP}
\setcounter{equation}{0}

The main result is the reverse duality of the non-conservative open ASEP
and the conservative shock ASEP with reflecting boundaries. The duality 
function is given by the Bernoulli shock measures of Definition 
\ref{Def:shockmeas} 
with $N$ shocks and a specific choice of shock densities determined by
condition \ref{bNshockcond2}. Since the main ideas in the proof for 
general $N$ appear already in the simpler case of a single shock, we 
separate the single-shock case (which also allows for more explicit 
information on the time evolution of the shock measure) from the general 
case. Before stating these results we assert reversibility of the shock 
exclusion process.

\begin{prop}
\label{Prop:revshockproc}
The $N$-particle shock exclusion process \ref{Def:shockprocess} on 
$\Lambda_L$ with reflecting boundaries is reversible w.r.t. the 
unnormalized measure
\be 
\pi(\bfx) = \prod_{i=1}^N d_i^{2x_i}
\label{revshockmeasN}
\ee
where
\be 
d_i := \sqrt{\frac{d^{r}_i}{d^{\ell}_i}}
\ee
is the hopping asymmetry of particle $i$.
\end{prop}

Now we are in a position to present the main results.

\begin{theo}
\label{Theo:TRD1shock} 
Let $W$ be the intensity matrix \eref{MASEP}
of the open ASEP with boundary rates \eref{lbr} and \eref{rbr} and 
for parameters $\rho_0$ and $\rho_1$ let
$Q$ be the intensity matrix of a simple biased random walk with generator
\eref{RWgen}.
Further, let $\mu^x$ be the shock 
measure \eref{shockmeas} with left bulk density $\rho_0=\rho_-$ and shock density 
\bel{rstar1}
\rho_1^{\star} = \frac{\alpha}{\alpha+\gamma}.
\ee 
The intensity matrices $W$ and $Q$ satisfy the reverse-duality relation
\be
\label{TRD1shock}
R W = Q^T R 
\ee
w.r.t. the duality matrix $R$ with matrix elements
$R_{x,\bfeta} = d_1^{2x} \mu^x(\bfeta)$ if and only if the following three 
conditions are satisfied:\\
(i) The bulk shock stability condition
\eref{shockcondi} is satisfied for $i=N=1$,\\ 
(ii) The boundary rates are on the manifold $\mathcal{B}^1_1$ .
\end{theo}

\begin{theo}
\label{Theo:Evolution1shock}
Denote by $\mu^x_t$
the distribution at time $t$ of the open ASEP, starting at $t=0$ from a 
shock measure $\mu^x_0=\mu^x$ and let Conditions (i) - (iii) of Theorem \ref{Theo:TRD1shock} be satisfied. Then, for any $x\in\Lambda_{L}$ 
\bel{Evolution1shock}
\mu^x(t) = \sum_{y = L_-}^{L_+}
P(y,t|x,0) \, \mu^y(0)
\ee
where 
\bel{condprobrw} 
P(y,t|x,0) =  \frac{d_1^2-1}{d_1^{2L}-1} d_1^{2(y-L_-)} 
+ \frac{2}{L} \sum_{p=1}^{L-1} d_1^{y-x} \psi_p(x) \psi_p(y)
 \frac{w }{\epsilon_p} \rme^{-\epsilon_p t}
\ee 
with
\bea
\psi_p(y) & := & d_1 \sin{\left(\frac{\pi p}{L}(y+1-L_-)\right)} 
- \sin{\left(\frac{\pi p}{L}(y-L_-)\right)} \\
\epsilon_p & = &  w \left[d_{1} + d^{-1}_{1} - 2 
\cos{\left( \frac{\pi p}{L} \right)}\right]   
\label{epsrw}
\eea
is the transition probability of the biased random walk generated by the intensity matrix $Q$ and starting at time $t=0$ from $x$.
The limit 
$\mu^\ast_1 := \lim_{t\to\infty} \mu^x(t)$ is the unique invariant 
measure
of the open ASEP with boundary parameters as specified by Conditions (ii) and (iii) and is given by the
convex combination
\bel{invmeas1}
\mu^\ast_1 =
\frac{d_1^2-1}{d_1^{2L}-1} \sum_{y=L_-}^{L_+} d_1^{2(y-L_-)}  \mu^y.
\ee
of shock measures $\mu^y$.
\end{theo}

\begin{rem}
For asymmetric exclusion processes with open boundaries
similar random walk dynamics were found earlier for a single anti-shock with 
bulk densities related by the inverse stability condition $z_1 = q^{-2} z_0$ 
under conditioning the process on a special atypical value of the 
time-integrated particle current \cite{Beli13b}. Boundary parameters in that 
work corresponded to $\omega_+=\omega_-=0$. In \cite{Jafa09} random 
walk dynamics of a single shock were proved for a discrete-time totally 
asymmetric simple exclusion process with deterministic sublattice update and 
stochastic open boundaries. In this model there is no constraint on the boundary parameters.
\end{rem}

\begin{rem}
According to condition (iii) in Theorem \ref{Theo:TRD1shock} the
invariant measure \eref{invmeas1} can be expressed by the two-dimensional 
representation of the stationary matrix product algebra characterized 
in the Appendices of \cite{Mall97} which exists for any $L\geq 2$ and which 
is analyzed further in \cite{Bryc19}. The homogeneous
Bernoulli measure with densities $\rho_0$ and $\rho_1$ generically
contribute to the two-dimensional MPM \cite{Jafa07} but are
absent in the convex combination \eref{invmeas1}. 
It would be interesting to explore whether two-dimensional
time-dependent matrices can capture the dynamics of some subspace
as has been proved for the SSEP with open boundaries \cite{Schu01}
and also for a non-conservative annihilating random walks with pair deposition \cite{Schu95,Schu20}.
\end{rem}

\begin{rem}
The spectrum of the generator given by the eigenvalues \eref{epsrw}
yields a subset of eigenvalues of the generator of the open ASEP
and is in agreement with the picture of spectral properties arising from a shock random walk explored in \cite{Dudz00,Sant02,deGi06}.
\end{rem}

As pointed out above, Theorem \ref{Theo:TRD1shock} is a special case of the following result for duality functions given by Bernoulli shock measures
with $N>1$ shocks.
\begin{theo}
\label{Theo:TRDNshock} 
Let $W$ be the intensity matrix \eref{MASEP}
of the open ASEP with boundary rates \eref{lbr} and \eref{rbr} and 
for parameters $\rho_0,\dots,\rho_N$ let
$Q$ be the intensity matrix of the $N$-particle shock exclusion process 
of Definition \ref{Def:shockprocess} and reversible measure \eref{revshockmeasN}.
Further, let $\mu^{\bfx}$ be the shock 
measure \eref{shockmeas} with left boundary density $\rho_0=\rho_-$ and shock fugacities
\bel{1shockcond}
z^{\star}_{i} = \frac{\alpha}{\gamma} q^{2(i-1)}
\ee 
for $1\leq i \leq N$.
The intensity matrices $W$ and $Q$ satisfy the reverse-duality relation
\be
\label{TRDNshock}
R W = Q^T R 
\ee
w.r.t. the duality matrix $R$ with matrix elements
$R_{x,\bfeta} = \pi(\bfx) \mu^{\bfx}_{\bfeta}$ if and only if the following three 
conditions are satisfied:\\
(i) The bulk shock stability condition
\eref{shockcondi} is satisfied for all $i\in\{1,\dots,N\}$,\\ 
(ii) The boundary rates are on the manifold $\mathcal{B}^1_N$ 
for $1\leq N \leq L$.
\end{theo}

\begin{rem}
The conservative reflective boundaries of the reverse dual are in contrast 
to the conventional duality for the open symmetric simple exclusion process 
(SSEP) which is dual to the SSEP with nonconservative absorbing boundaries 
noticed in \cite{Spoh83} and later fully developed in \cite{Cari13,Fras20}, 
including the generalization to symmetric partial exclusion \cite{Schu94}. 
An open question is whether a similar conservative reverse duality can be 
constructed for the conservative asymmetric partial exclusion process studied 
in \cite{Cari16} when non-conservative open boundaries are added.
\end{rem}

Since the shock exclusion process with reflecting boundaries is reversible for any $N$ also Theorem \ref{Theo:Evolution1shock} 
has a natural generalization for $N>1$.

\begin{theo}
\label{Theo:Nshockevolution}
Let $\mu^{\bfx}(t)$ denote
the distribution at time $t$ of the ASEP of Theorem \ref{Theo:TRDNshock}, starting from an
$N$-shock measure $\mu^\bfx=\mu^\bfx(0)$. Then
\bel{shockevolutionN}
\mu^{\bfx}(t) = \sum_{\bfy} 
P(\bfy,t|\bfx,0) \, \mu^{\bfy}
\ee
where $P(\bfy,t|\bfx,0)$ is the transition probability 
of the shock exclusion process \ref{Def:shockprocess}. The limit
$\mu^\ast_N := \lim_{t\to\infty}\mu^{\bfx}(t)$ is the unique invariant 
measure
of the open ASEP with boundary parameters as specified by Conditions (ii) and (iii) and is given by the unnormalized
convex combination
\bel{invmeasN}
\tilde{\mu}^\ast_N
= \frac{j_{N}^{2y_{N}}}{j_{1}^{2y_{1}}} \sum_{\bfy}  \prod_{i=1}^{N-1} j_{i}^{2(y_i-y_{i+1})} \mu^{\bfy}
\ee
of shock measures with $N$ shocks.
\end{theo}

\begin{rem}
The spectrum of the generator of the $N$-particle shock exclusion process 
can be calculated from Bethe ansatz \cite{Nepo03b,Nepo04,deGi05,deGi06,Simo09}. However, there is no simple expression in closed form. For results on the eigenfunctions see \cite{Cram10}. Further progress along these lines has been made recently
for the SSEP with open boundaries by using duality \cite{Fras20}.
It would be most interesting to generalize this approach to the open ASEP to cover the
full manifold of boundary parameters.
\end{rem}

\section{Proofs}
\setcounter{equation}{0}

We prove reverse duality and its consequences for the time evolution of
shock measures by studying the action of the quantum 
Hamiltonian on the vector representation of the shock measures
which are shown to form an invariant subspace.
The coefficients appearing in these computations yield the reverse duality
with the shock exclusion process.
Sufficient conditions for a similar property were
conjectured earlier in \cite{Kreb03} for second-class particles
at the shock positions and in \cite{Jafa07} for a single shock.
Here we obtain rigorously the necessary and sufficient
conditions on the parameters of the ASEP for the reverse duality for a 
specific choice of shock densities for $N\geq 1$ shocks.. 

\subsection{Preliminaries}

Everywhere below a vector denoted
by $\ket{abc\dots}$ is defined to be the Kronecker product 
$\ket{a}\otimes\ket{b}\otimes\ket{c}\otimes\dots$ of two-dimensional
vectors $\ket{a},\ket{b},\ket{c}\dots \in \R^2$ with components $(a_0,a_1), (b_0,b_1), (c_0,c_1), \dots$. For vectors with components of the form $a_0 = 1-\rho$, $a_1 = \rho$ and $0\leq\rho\leq 1$ we use greek letters and write $\ket{\rho}$
and call $\rho$ a density.
In slight abuse of wording we call for any vector $\ket{a}$ with $a_0\neq 0$ the ratio
\be
z(a) := \frac{a_1}{a_0}
\label{zgendef}
\ee
the {\it fugacity}. The functions
\be 
\epsilon_-(\cdot) := \alpha - \gamma z(\cdot), \quad
\epsilon_+(\cdot) := \delta - \beta z(\cdot)
\label{epspmdef}
\ee
play a role as eigenvalues of the boundary operators $h^\pm$.

For two vectors $\ket{a}$, $\ket{b}$ we introduce the determinant
\be
\Delta(a,b) := a_0 b_1 - a_1 b_0
\label{detdef}
\ee
and recall that linear dependence of $\ket{a}$ and $\ket{b}$
is equivalent to $\Delta(a,b) = 0$. For $z(b) = q^2 z(a)$ 
one has
\be
\Delta(a,b) = (r-\ell) \frac{a_0 b_1}{r}  = (r-\ell) \frac{a_1 b_0}{\ell}.
\label{Deltaq2}
\ee
For linearly independent vectors
we also define the functions
\bea
& & d(a,b) := (r-\ell)\frac{a_1a_0}{\Delta(a,b)}, \quad
\tilde{d}(a,b) := (r-\ell)\frac{a_1b_0}{\Delta(a,b)}.
\label{dbulkdef} \\
& & d_\pm(a,b) := \frac{\epsilon_\pm(a) a_0 (a_0+a_1)}{\Delta(a,b)}, 
\quad \tilde{d}_\pm(a,b) := \frac{\epsilon_\pm(a) a_0 (b_0+b_1)}{\Delta(a,b)} .
\label{dpmdef}
\eea
which play a role for the transition rates of the dual process.
We recall that in all assertions and in all the proofs below it is tacitly assumed 
that all boundary rates and bulk jump rates of the open ASEP
are strictly positive and that $q\neq 1$.

\begin{lmm}[Eigenvectors]
\label{Lem:eigen}
(i) For any two vectors $\ket{a}$, $\ket{b}$ the bulk eigenvalue property
\bel{hib}
h \ket{ab} = 0
\ee
is satisfied if and only if $\Delta(a,b)=0$.\\
(ii) For vectors $\ket{a^\pm}$ with components such that
 $a^\pm_0 a^\pm_1 \neq 0$ and $a^\pm_0 + a^\pm_1 \neq 0$
the two boundary eigenvalue equations 
\be 
h^\pm \ket{a^\pm} = \epsilon_\pm(a^\pm) \ket{a^\pm} 
\label{eigen}
\ee 
have a unique solution 
with strictly positive fugacities given by
\be 
z(a^-) = \kappa^{-1}(\alpha,\gamma),
\quad z(a^+) = \kappa_+(\beta,\delta)
\label{eigenz}
\ee 
with $\kappa_+(\cdot,\cdot)$ defined in \eref{kappadef}. The eigenvalues are
\bea
\epsilon_-(\cdot) & = & \alpha - \gamma z(\cdot) \, = \, (r-\ell) \frac{z(\cdot)}{1+z(\cdot)}
\label{epsm} \\
\epsilon_+(\cdot) & = & \delta - \beta z(\cdot) \, = \, - (r-\ell)
\frac{z(\cdot)}{1+z(\cdot)}.
\label{epsp}
\eea
For $a^\pm_0 a^\pm_1=0$ or $a^\pm_0 + a^\pm_1 = 0$ the eigenvalue equation \eref{eigen} has no solution for strictly positive boundary rates
and $r\neq \ell$.\\
\end{lmm}

\Proof Consider two vectors $\ket{a}$ with components $a_0,a_1$
and $\ket{b}$ with components $b_0,b_1$ and determinant 
$\Delta_{ab}$ defined in \eref{detdef}.\\
(i) For the product vector one gets from the definition \eref{hbulk}
\bea
h \ket{ab}
& = &
\left(\ba{cccc}
0 & 0 & 0 & 0 \\
0 & r & - r & 0 \\
0 & - \ell & \ell & 0 \\
0 & 0 & 0 & 0 
\ea \right) 
\left(\ba{c}
a_0 b_0\\
a_0 b_1 \\
a_1 b_0 \\
a_1 b_1
\ea \right) =
\Delta_{ab} \left(\ba{c}
0\\
r \\
- \ell \\
0
\ea \right) 
\label{hibcomp}
\eea
which proves part (i) of the Lemma.\\
(ii) To make notation less heavy we drop the superscripts $\pm$
in the vectors $\ket{a^\pm}$ in the following
computations. From the definition \eref{hbm} of the boundary matrix
one obtains
\be
h^- \ket{a} = \left(\ba{c} 
\alpha a_0 - \gamma a_1 \\
- \alpha a_0 + \gamma a_1 + (r-\ell) a_1
\ea\right) 
\ee
from which one reads off immediately that for 
strictly positive rates $\alpha,\gamma$ there exists no eigenvector with 
$a_0=0$ or $a_1=0$, i.e., for
$a_0 a_1 =0$. To proceed with the proof of the eigenvalue properties 
for a vector $\ket{a}$ such that $a_0 a_1 \neq 0$ we 
define
\be 
g_\pm(a) := \mp \frac{r-\ell}{\epsilon_-(a)} - \frac{1}{z(a)}
\ee
and 
note that for any linearly independent vector 
$\ket{b}$ the equality
\bea 
h^- \ket{a} 
& = & \epsilon_-(a) \left(\ba{c} 
a_0  \\
g_-(a) a_1
\ea\right) \nonumber \\
& = & \left(\tilde{d}_-(a,b) - \tilde{d}(a,b)\right) \ket{a} - \left(d_-(a,b) - d(a,b)\right)  \ket{b}
\label{hmcomp}
\eea
holds without further constraints on the vectors $\ket{a}$ and $\ket{b}$.

If an eigenvector exists (which requires $d_-(a,b) = d(a,b)$
or equivalently $g_-(a) =1$)
then the eigenvalue is of the form $\epsilon_-(a)$ as defined in 
\eref{epspmdef} and one deduces that there exists no eigenvector with $a_1 = - a_0$ 
since then $z=-1$ and therefore $g(a) \neq 1$ for $r\neq\ell$.
From $g_-(a) =1$ one also obtains the second equality in \eref{epsm} 
and together with the first equality one arrives at the quadratic equation
\be  
z^2 + \frac{\gamma - \alpha + r - \ell}{\gamma} z - \frac{\alpha}{\gamma} = 0
\ee 
for the fugacity $z\equiv z(a)$.
The positive root is given by $z = \kappa^{-1}(\alpha,\gamma)$
which yields (up to an arbitrary multiplicative factor) the components 
of the eigenvector $\ket{a}$ in terms  of the boundary rates as in 
\eref{eigenz}.

For the right boundary matrix one gets 
\bea
h^+ \ket{a} & = & \left(\ba{c} 
\delta a_0 - \beta a_1 \\
- \delta a_0 + \beta a_1 - (r-\ell) a_1
\ea\right) \nonumber \\
& = & \epsilon_+(a) \left(\ba{c} 
a_0 \\
g_+(a) a_1
\ea\right) \nonumber \\
& = & \left(\tilde{d}_+(a,b) + \tilde{d}(a,b)\right) \ket{a} - \left(d_+(a,b) + d(a,b)\right)  \ket{b}.
\label{hpcomp}
\eea
with a linearly independent vector $\ket{b}$.
Eigenvectors with $a_0a_1=0$ or
$a_0+a_1=0$ cannot exist as can be seen by arguments analogous to the
case treated above. For an eigenvector one immediately reads off the first equality in \eref{epsp}. Moreover, 
an eigenvector satisfies $d_+(a,b) = - d(a,b)$ which gives
the second equality in \eref{epsp} and
together with the first equality in \eref{epsp} the quadratic equation
\be  
z^2 - \frac{\delta - \beta + r - \ell}{\beta} z - \frac{\delta}{\beta} = 0 
\ee 
for the fugacity. The
positive root is given by $z = \kappa_+(\beta,\delta)$ 
which yields  \eref{eigenz} for the right boundary.
\qed

\begin{cor}
\label{Cor:Nshockcond}
A probability vector $\ket{\pi}$ with factorized boundary marginals such that $\ket{\pi} = \ket{\rho_-} \otimes \ket{\tilde{\pi}}\otimes
\ket{\rho_+}$ and with
$z_+ = q^{2N} z_-$
is an eigenvector of the boundary operator
$h^-_{L_-} + h^+_{L_+}$ if and only if \eref{bNshockcond} holds.
\end{cor}

\begin{lmm}[Projection]
\label{Lem:hbulk} 
(i) For vectors $\ket{a^{(i)}}$,  $\ket{\tilde{a}^{(i)}}$, $i\in\{1,2\}$, 
such that $a^{(i)}_0  a^{(i)}_1 \neq 0$ and $ \tilde{a}^{(i)}_0  
\tilde{a}^{(i)}_1 \neq 0$
and a pair of linearly independent vectors
($\ket{b^{(i)}}$, $\ket{\tilde{b}^{(i)}}$) such that
$\ket{b^{(1)}}$ is linearly independent of $\ket{\tilde{a}^{(1)}}$ 
and $\ket{b^{(2)}}$ is linearly independent of
$\ket{\tilde{a}^{(2)}}$ the two-site jump generator 
$h$ has the three-dimensional projection properties 
\bea
(A) \qquad h \ket{a^{(1)}b^{(1)}} 
& = & u^{(1)} \ket{a^{(1)}b^{(1)}} 
- v^{(1)} \ket{b^{(1)}\tilde{a}^{(1)}} 
- w^{(1)} \ket{a^{(1)}\tilde{b}^{(1)}} 
\label{hbulk1} 
\eea 
with coefficients
\be
u^{(1)} = \ell - \tilde{d}(b^{(1)},\tilde{b}^{(1)}), \quad
v^{(1)} = d(a^{(1)},\tilde{a}^{(1)}), \quad
w^{(1)} = - d(b^{(1)},\tilde{b}^{(1)})
\label{hbulk1coeff} 
\ee
and
\bea
(B) \qquad h \ket{b^{(2)}\tilde{a}^{(2)}} 
& = & u^{(2)} \ket{b^{(2)}\tilde{a}^{(2)}}  - v^{(2)}  \ket{a^{(2)}b^{(2)}} 
- w^{(2)} \ket{\tilde{b}^{(2)}\tilde{a}^{(2)}} 
\label{hbulk3}
\eea
with coefficients
\be
u^{(2)} = r + \tilde{d}(b^{(2)}, \tilde{b}^{(2)}), \quad
v^{(2)} = - d(\tilde{a}^{(2)}, a^{(2)}), \quad
w^{(2)} = d(b^{(2)}, \tilde{b}^{(2)})
\label{hbulk3coeff} 
\ee
if and only if the condition
\be 
\frac{z(\tilde{a}^{(i)})}{z(a^{(i)})} = q^2, \quad 
\label{hab2} 
\ee
on the 
vectors $\ket{a^{(i)}}$,  $\ket{\tilde{a}^{(i)}}$ are satisfied. \\
(ii) For a pair of linearly independent 
vectors ($\ket{a}$,  $\ket{\tilde{a}}$) such that $a_0  a_1 \neq 0$ and $\tilde{a}_0  
\tilde{a}_1 \neq 0$
and vectors
$\ket{\tilde{c}}$ linearly independent of $\ket{a}$, $\ket{c}$
linearly independent of $\ket{\tilde{a}}$ the two-site jump generator $h$ has the three-dimensional
projection property 
\bea
(C) \qquad h \ket{a\tilde{a}}
& = & 
u \ket{a\tilde{a}} - w \ket{\tilde{c}\tilde{a}} - v \ket{a c}
\label{hbulk2} 
\eea
with coefficients
\be 
u = \tilde{d}(a,\tilde{c}) - \tilde{d}(\tilde{a},c), \quad
v = - d(\tilde{a},c), \quad w = d(a,\tilde{c})
\label{hbulk2coeff} 
\ee
if and only if  the condition \eref{hab2} is satisfied.\\
(iii) If $a_0  a_1 = 0$ or $ \tilde{a}_0  
\tilde{a}_1 = 0$ the projection equations \eref{hbulk1}, \eref{hbulk3}, 
and \eref{hbulk2} have no solution for strictly positive boundary
rates and $r\neq \ell$.
\end{lmm}

\Proof We focus on the projection property (A) and drop the superscript $(1)$
to make notation less heavy. The proof of the projection properties (B) and (C) is analogous. In simplied notation \eref{hbulk1}
reads $h \ket{ab} = u \ket{ab} - v \ket{b\tilde{a}} - w \ket{a\tilde{b}}$
from which using \eref{hibcomp} one obtains are the four linear equations
\bea 
0 
& = & u a_0 b_0 
- v b_0 \tilde{a}_0  
- w a_0 \tilde{b}_0 \\
r (a_0 b_1 - a_1 b_0)
& = & u a_0 b_1 
- v b_0 \tilde{a}_1  
- w a_0 \tilde{b}_1 \\
- \ell (a_0 b_1 - a_1 b_0)
& = & u a_1 b_0 
- v b_1 \tilde{a}_0  
- w a_1 \tilde{b}_0 \\
0 
& = & u a_1 b_1 
- v b_1 \tilde{a}_1  
- w a_1 \tilde{b}_1 
\eea
for the three coefficients $u,v,w$.
(a) Consider $a_0 = 0$ which implies by linear independence
$b_0 \neq 0$. Then the first two equations yield $v=r=0$
in contradiction to $r>0$. Hence there exists no solution for
$a_0 = 0$. \\
(b) Consider $a_0 \neq 0$ and $\tilde{a}_0 = 0$. The 
first and third equation reduce to
\bea 
0 
& = & u b_0 
- w \tilde{b}_0 \\
- \ell (a_0 b_1 - a_1 b_0)
& = & a_1  \left(u b_0 
- w \tilde{b}_0 \right)
\eea
which due to linear independence yields $\ell = 0$, in contradiction
to the assumption $\ell >0$. Hence there exists no solution for
$\tilde{a}_0 = 0$. The proof that there exists no solution
for $a_1 = 0$ or $\tilde{a}_1 = 0$ is analogous.\\
(c) To proceed with the generic case we introduce the fugacities
\be 
\tilde{z} := \frac{\tilde{a}_1}{\tilde{a}_0}, \quad
z := \frac{a_1}{a_0} 
\ee
and consider next $b_0 = 0$. Using linear independence the four equations yield immediately the unique solution
\bea 
u = r , \quad
v = \ell \frac{a_0}{\tilde{a}_0}, \quad
w = 0
\eea
if and only if $\tilde{z} = q^2 z$
in agreement with the assertion of the theorem.
For $b_0 \neq 0$
it is useful to define
\be
\tilde{v} := v \frac{\tilde{a}_0}{a_0}, \quad 
\tilde{w} := w \frac{\tilde{b}_0}{b_0} .
\ee
The four equations turn into
\bea 
u & = & \ell  + \tilde{w} \\
\tilde{v} & = & \ell 
\eea
and
\bea
r (a_0 b_1 - a_1 b_0)
& = & a_0 \left[(\ell  + \tilde{w}) b_1 
- \tilde{w} \frac{b_0}{\tilde{b}_0} \tilde{b}_1\right]
- \ell \frac{a_0}{\tilde{a}_0} b_0 \tilde{a}_1  \\
0 & = & a_1 \left[(\ell  + \tilde{w}) b_1 
- \tilde{w} \frac{b_0}{\tilde{b}_0} \tilde{b}_1\right]
- \ell \frac{a_0}{\tilde{a}_0} b_1 \tilde{a}_1  .
\eea
Inserting the second into the first equation then yields
\bea
\tilde{w} & = & \frac{(r - \ell) b_1 \tilde{b}_0}{b_1 \tilde{b}_0 - b_0 
\tilde{b}_1} \, = \, - \frac{(r - \ell) b_1 \tilde{b}_0}{\Delta_{b\tilde{b}}}
\eea
and this solution exists
if and only if $\tilde{z} = q^2 z$ as the second equation shows. 
Thus with \eref{detdef} one arrives at
\bea
u & = & \ell - (r - \ell) \frac{b_1 \tilde{b}_0}{\Delta_{b\tilde{b}}}  \\
v & = & \ell \frac{a_0}{\tilde{a}_0} \, = \, 
(r-\ell) \frac{a_0 a_1}{\Delta_{a\tilde{a}}} \\
w & = & - (r - \ell) \frac{b_0 b_1}{\Delta_{b\tilde{b}}}
\eea
and with the functions $d(\cdot,\cdot)$ and $\tilde{d}(\cdot,\cdot)$ defined 
in \eref{dbulkdef} one obtains the coefficients \eref{hbulk1coeff}.
\qed

\begin{cor} 
\label{Cor:shockjumprates}
For shocks satisfying the shock stability criterion \eref{shockcondi} 
the shock 
jump rates $d_{ir,\ell}$ satisfy
the identities
\bea 
d(\rho_{i-1},\rho_{i}) & = & \ell \frac{1-\rho_{i-1}}{1-\rho_{i}} 
\, = \, r  \frac{\rho_{i-1}}{\rho_{i}} 
\, = \, \ell (1-\rho_{i-1}) + r \rho_{i-1} \, = \, d^{\ell}_{i} 
\label{vi} \\
- d(\rho_{i},\rho_{i-1}) & = & \ell \frac{\rho_{i}}{\rho_{i-1}} 
\, = \,
r \frac{1-\rho_{i}}{1-\rho_{i-1}} 
\, = \, r (1-\rho_{i}) + \ell \rho_i \, = \, d^{r}_{i} 
\label{wi} \\
& & \epsilon_+(\rho_i) \, = \, d^{r}_{i} - r 
\label{epspi} \\
& & \epsilon_-(\rho_{i-1}) \, = \, d^{\ell}_{i} - \ell .
\label{epsmi}
\eea
which also imply $d^{r}_{i} d^{\ell}_{i} = r\ell$ independently of $i$.
\end{cor}

\subsection{Proof of Proposition \ref{Prop:revshockproc}}

The proof of reversibility of the shock exclusion process) is essentially trivial. According to Definition \ref{Def:shockprocess} the intensity matrix of the $N$-particle
shock exclusion process is given by the matrix elements
\be 
Q_{\bfx\bfx'} = \sum_{i=1}^{N} \left[w^\ell_{i}(\bfx) (\delta_{\bfx',\bfx^{i-}} - \delta_{\bfx',\bfx}) + w^r_{i}(\bfx) (\delta_{\bfx',\bfx^{i+}} - \delta_{\bfx',\bfx})\right].
\label{shockasepQ}
\ee
Using $\delta_{\bfy^{i\pm},\bfx} = \delta_{\bfy,\bfx^{i\mp}}$
one obtains
\bea
d_i^2 w^{\ell}_{i}(\bfx') \delta_{\bfx,(\bfx')^{i-}}
& = & w^{r}_{i}(\bfx) \delta_{\bfx',\bfx^{i+}} \\
d_i^{-2} w^{r}_{i}(\bfx') \delta_{\bfx,(\bfx')^{i+}}
& = & w^{\ell}_{i}(\bfx) \delta_{\bfx',\bfx^{i-}}
\eea
and $\pi(\bfx^{i\pm})/\pi(\bfx) = d_i^{\pm 2}$ yields
$Q_{\bfy\bfx} = Q_{\bfx\bfy} \pi(\bfx)/\pi(\bfy)$ which is the 
definition of reversibility, viz.,
$Q^T = \hat{\pi} Q \hat{\pi}^{-1}$,
in terms of the matrix elements.
\qed

\subsection{Reverse duality involving one shock}

The modified duality matrix $S$ with elements $S_{x,\bfeta} = 
\inprod{\bfeta}{\mu^x} = \mu^x(\bfeta)$ is the matrix with the column 
vectors $\ket{\mu^{x}}$ as rows. 
For one shock the vector representation \eref{shockmeasvec} of the shock measure of Theorem \ref{Theo:measure1} is given by
\bel{shockmeasvec1}
\ket{\mu^x} = \ket{\rho_{0}}^{\otimes (x - L_-)} \otimes \ket{\rho^\star} \otimes \ket{\rho_{1}}^{\otimes (L_+ - x)}, \quad
L_- \leq x \leq L_+
\ee 
with $\rho^\star$ given by \eref{rstar1}. 

\subsubsection{Proof of Theorem \ref{Theo:TRD1shock}}

One has to find the necessary
and sufficient conditions under which the shock measures
$\ket{\mu^{x}}$ form an invariant subspace under the action of $H$, i.e.,
\be 
H \ket{\mu^{x}} = \sum_{y=L_-}^{L_+} G_{yx} \ket{\mu^{y}}
\label{Hmu1}
\ee
with a matrix $G$ and then show that $G=-Q^T$ for
the intensity matrix $Q$ of the biased random walk defined in
Remark \ref{Rem:RW}. To do so, we first employ the projection 
lemma \ref{Lem:hbulk} together with the eigenvector
properties detailed in Lemma \ref{Lem:eigen} and its Corollary 
\ref{Cor:Nshockcond} to prove \eref{Hmu1}
if and only if \eref{bNshockcond}, \eref{bNshockcond2} and \eref{Nshockcond} hold with $N=m+1=1$. 

\paragraph{Shock in the range $L_- < x < L_+$:} From the product
form of the shock measure, the eigenvalue properties \eref{hib}
and \eref{eigen}
in Lemma \ref{Lem:eigen} and Corollary \ref{Cor:Nshockcond} one concludes
\bea 
H \ket{\mu^{x}} & = & (\tilde{h}^{-}_{L_-} + \tilde{h}_{x-1,x} +\tilde{h}_{x,x+1}
+ \tilde{h}^{+}_{L_+}) \ket{\mu^{x}} \nonumber \\
& = & \ket{\rho_{0}}^{\otimes (x - L_- -1)}\otimes
[(\tilde{h}_{12} + \tilde{h}_{23}) \ket{\rho_{0}} \otimes \ket{\rho_1^{\star}}\otimes
\ket{\rho_{1}}] \otimes \ket{\rho_{1}}^{\otimes (L_+ - x -1)}
\nonumber \\
& & + (\epsilon_-(\rho_0) + \epsilon_+(\rho_1)) \ket{\mu^{x}}.
\eea
Taking $a^{(1)}=a^{(2)}=\rho_0$, $\tilde{a}^{(1)} = \tilde{a}^{(2)} = \rho_1$,
$b^{(1)} = b^{(2)} = \rho_1^\star$ in \eref{hbulk1} and 
\eref{hbulk3} of Lemma \ref{Lem:hbulk}
yields after taking Kronecker products on the right in \eref{hbulk1}
and on the left in \eref{hbulk3}
\bea
\tilde{h}_{12} \ket{\rho_0\rho_1^\star\rho_1} 
& = & \left(\ell - \tilde{d}(\rho_1^\star,\tilde{b}^{(1)})\right) \ket{\rho_0\rho_1^\star\rho_1} \nonumber \\
& & - d(\rho_0,\rho_1) \ket{\rho_1^\star\rho_1\rho_1} 
+ d(\rho_1^\star,\tilde{b}^{(1)}) \ket{\rho_0\tilde{b}^{(1)}\rho_1} 
\nonumber \\
\tilde{h}_{23} \ket{\rho_0\rho_1^\star\rho_1} 
& = & \left(r + \tilde{d}(\rho_1^\star, \tilde{b}^{(2)}) \right) \ket{\rho_0\rho_1^\star\rho_1} \nonumber \\
& & + d(\rho_1, \rho_0) \ket{\rho_0\rho_0\rho_1^\star} 
- d(\rho_1^\star, \tilde{b}^{(2)}) \ket{\rho_0\tilde{b}^{(2)}\rho_1} .
\eea
Since $\tilde{b}^{(1)}$ and $\tilde{b}^{(2)}$ are arbitrary (except for 
the immaterial inequality with $\rho_1^{\star}$) we can without loss of 
generality take $\tilde{b}^{(1)} = \tilde{b}^{(2)} \neq \rho^{\star}$. 
Thus the projection lemma asserts that for $q>1$, $\rho_1^{\star}\in[0,1]$ 
and any pair of densities $\rho_0, \rho_{1} \in (0,1)$ such that 
$0<\rho_{0}<\rho_{1}<1$ and $z_{1} = q^2 z_{0}$ the three-site jump 
generator $\tilde{h}_{12} + \tilde{h}_{23}$ projects the vector $\ket{\rho_0\rho_1^{\star}\rho_1} := \ket{\rho_0} \otimes \ket{\rho_1^{\star}} \otimes \ket{\rho_{1}}$
on the three-dimensional subspace of spanned by $\ket{\rho_0\rho_1^{\star}\rho_1}$,
$\ket{\rho_1^{\star}\rho_1\rho_1} := \ket{\rho_1^{\star}} \otimes \ket{\rho_1}^{\otimes 2}$, and $\ket{\rho_0\rho_0\rho_1^{\star}} := \ket{\rho_0}^{\otimes 2} \otimes \ket{\rho_1^{\star}}$ as
\be 
(\tilde{h}_{12} + \tilde{h}_{23})\ket{\rho_0\rho_1^{\star}\rho_1} =
(r+\ell) \ket{\rho_0\rho_1^{\star}\rho_1} 
-  v_{01} \ket{\rho_1^{\star}\rho_1\rho_1}
- w_{01} \ket{\rho_0\rho_0\rho_1^{\star}}
\label{hbulkclosure}
\ee
with constants
\bel{uvw} 
u_{01} = r+\ell, \quad w_{01} = - d(\rho_{1},\rho_{0}), \quad 
v_{01} = d(\rho_{0},\rho_{1}) 
\ee
that do not depend on $\rho_1^{\star}$. The identities \eref{vi} -
\eref{epsmi} then yield for arbitrary shock density $\rho_1^{\star}$
the bulk closure relation
\bea 
H \ket{\mu^{x}} 
& = & [r + \ell + \epsilon_-(\rho_0) + \epsilon_+(\rho_1)] \ket{\mu^{x}}
 - d^{\ell}_{1} \ket{\mu^{x-1}} - d^{r}_{1} \ket{\mu^{x+1}} \nonumber \\
& = & (d^{\ell}_{1} + d^{r}_{1}) \ket{\mu^{x}}
 - d^{\ell}_{1} \ket{\mu^{x-1}} - d^{r}_{1} \ket{\mu^{x+1}}
\eea
if and only if the conditions \eref{bNshockcond} and \eref{Nshockcond}
are satisfied with $N=1$.

\paragraph{Shock at the left boundary:} 
From the eigenvalue properties in Lemma \ref{Lem:eigen} one gets
\bea 
H \ket{\mu^{L_-}} 
& = & (h^-_{L_-} + \tilde{h}_{L_-,L_-+1} + h^+_{L_+})
\ket{\mu^{L_-}} \nonumber \\
& = & [(h^-_{1} + h) \ket{\rho_1^{\star}}\otimes \ket{\rho_1}] \otimes
\ket{\rho_1}^{\otimes{L-2}}
+ \epsilon_+(\rho_1) \ket{\mu^{L_-}} .
\label{HLm} 
\eea
To handle the left boundary term in \eref{HLm} use case B in the 
projection lemma \ref{Lem:hbulk} by
choosing in \eref{hbulk3} $a^{(2)}=\rho_{0}$, $b^{(2)} = \rho_1^{\star}$, and 
$\tilde{a}^{(2)} = \rho_{1}$ and by choosing in \eref{hmcomp}
$a = \rho_1^{\star}$ taking the Kronecker product with $\ket{\rho_1}$ 
on the right. This yields
\bea 
h \ket{\rho_1^{\star}\rho_{1}} 
& = & \left(r + \tilde{d}(\rho_1^{\star}, \tilde{b}^{(2)}) \right) \ket{\rho_1^{\star}\rho_{1}} \nonumber \\
& & + d(\rho_{1}, \rho_{0}) \ket{\rho_{0}\rho_1^{\star}} 
- d(\rho_1^{\star}, \tilde{b}^{(2)}) \ket{\tilde{b}^{(2)}\rho_{1}} 
\eea
and
\bea
\tilde{h}_1 \ket{\rho_1^{\star}\rho_{1}}
& = & \left(\tilde{d}_-(\rho_1^{\star},b) - \tilde{d}(\rho_1^{\star},b)\right) \ket{\rho_1^{\star}\rho_{1}} - \left(d_-(\rho_1^{\star},b) - d(\rho_1^{\star},b)\right)  \ket{b\rho_{1}}
\eea
Taking without loss of generality $b=\tilde{b}^{(2)}\neq\rho_1^{\star}$ one finds
\bea 
(h + \tilde{h}_1) \ket{\rho_1^{\star}\rho_{1}} 
& = & \left(r + \tilde{d}_-(\rho_1^{\star}, b) \right) 
\ket{\rho_1^{\star}\rho_{1}} \nonumber \\
& & + d(\rho_{1}, \rho_{0}) \ket{\rho_{0}\rho_1^{\star}} 
- d_-(\rho_1^{\star}, b) \ket{b\rho_{1}} .
\eea
Observe next that for $\rho_1^\star$ as specified in the Theorem one has $\epsilon_-(\rho_1^\star)=0$ so that $\tilde{d}_-(\rho_1^{\star}, b) = d_-(\rho_1^{\star}, b) = 0$ 
for any vector $b\neq \rho_1^{\star}$. Therefore 
\bea 
(h + \tilde{h}_1) \ket{\rho_1^{\star}\rho_{1}} 
& = & r \ket{\rho_1^{\star}\rho_{1}} 
+ d(\rho_{1}, \rho_{0}) \ket{\rho_{0}\rho_1^{\star}} .
\eea
Using \eref{wi} and \eref{epspi} the conclusion from these computations 
is the left closure relation
\be 
H \ket{\mu^{L_-}} = [r+\epsilon_+(\rho_1)] \ket{\mu^{L_-}} 
+ d(\rho_{1}, \rho_{0}) \ket{\mu^{L_-+1}} 
=
d^{r}_{1} \left(\ket{\mu^{L_-}} - \ket{\mu^{L_-+1}}\right)
\ee 
if and only if $\rho_1^\star = \alpha/(\alpha+\gamma)$ as specified in the 
Theorem.

\paragraph{Shock at the right boundary:}

From Lemma \ref{Lem:eigen} one gets in a similar fashion
\bea 
H \ket{\mu^{L_+}} & = & (h^-_{L_-} + \tilde{h}_{L_+-1,L_+} + h^+_{L_+})
\ket{\mu^{L_+}} \nonumber \\
& = &  \epsilon_-(\rho_0) \ket{\mu^{L_+}} 
+ \ket{\rho_0}^{\otimes(L-2)} 
\otimes [(h + h^+_2)
\ket{\rho_0} \otimes \ket{\rho^{\star}}]
\label{HLp}
\eea
To treat \eref{HLp} we use the projection formula
\eref{hbulk1} in Lemma \ref{Lem:hbulk}
with the replacements $a^{(1)} \to \rho_0$, $b^{(1)} \to \rho_1^{\star}$, $\tilde{a}^{(1)} \to \rho_1$ and take 
$a = \rho_1^{\star}$ 
in \eref{hpcomp}. After left Kronecker multiplication
with $\ket{\rho_0}$ this yields
\bea 
h \ket{\rho_0\rho_1^{\star}} 
& = & \left(\ell - \tilde{d}(\rho_1^{\star},\tilde{b}^{(1)})\right) \ket{\rho_0\rho_1^{\star}} \nonumber \\
& & - d(\rho_0,\rho_1) \ket{\rho_1^{\star}\rho_1} 
+ d(\rho_1^{\star},\tilde{b}^{(1)}) \ket{\rho_0\tilde{b}^{(1)}} 
\eea
and
\be
h^+_2 \ket{\rho_0\rho_1^{\star}} = \left(\tilde{d}(\rho_1^{\star},b) + d_+(\rho_1^{\star},b)\right)  \ket{\rho_0\rho^{\star}} 
- \left(d(\rho_1^{\star},b) + d_+(\rho_1^{\star},b)\right) \ket{\rho_0b} 
\ee 
With $b=\tilde{b}^{(1)}\neq \rho_0$ one thus finds the projection
property
\bea 
(h + h^+_2)\ket{\rho_0\rho_1^{\star}} 
& = & \left(\ell - \tilde{d}(\rho_1^{\star},\tilde{b}^{(1)})
+\tilde{d}(\rho_1^{\star},b) + d_+(\rho_1^{\star},b)\right) \ket{\rho_0\rho_1^{\star}} \nonumber \\
& & - d(\rho_0,\rho_1) \ket{\rho_1^{\star}\rho_1} 
+ d(\rho_1^{\star},\tilde{b}^{(1)}) \ket{\rho_0\tilde{b}^{(1)}} 
\nonumber \\
& & 
- \left(d(\rho_1^{\star},b) + d_+(\rho_1^{\star},b)\right) \ket{\rho_0b} 
\eea
which with \eref{wi} and \eref{epsmi} and by choosing without 
loss of generality $b=\tilde{b}^{(1)}\neq\rho_1^{\star}$ leads to 
\bea 
(h + h^+_2)\ket{\rho_0\rho_1^{\star}} 
& = & \left(\ell + d_+(\rho_1^{\star},b)\right) \ket{\rho_0\rho_1^{\star}} \nonumber \\
& & - d(\rho_0,\rho_1) \ket{\rho_1^{\star}\rho_1} 
- d_+(\rho_1^{\star},b) \ket{\rho_0b} 
\eea
Since by assumption \ref{bNshockcond2} holds one has $\epsilon_+(\rho_1) = 0$ so that
$d_+(\rho_1^{\star},b) = 0$ for all $b\neq\rho_1^{\star}$.
With \eref{vi} we arrive at the right closure relation
\bea 
H \ket{\mu^{L_+}} 
& = &  d^{\ell}_{1}\left(\ket{\mu^{L_+}} - \ket{\mu^{L_+-1}}\right).
\eea
Conversely, if \ref{bNshockcond2} does not hold then closure is violated.

It is thus proved that the set of vectors $\ket{\mu^x}$, $x\in\Lambda_L$
forms an invariant subspace under the action of $H$ if and only of the
conditions (i) - (iii) hold.
This allows for writing \eref{Hmu1}
with coefficients 
\be 
G_{yx}
= \left\{ \ba{ll}
d^{r}_{1} \left(\delta_{y,L_-} - \delta_{y,L_-+1}\right) & x = L_- \\
(d^{\ell}_{1} + d^{r}_{1}) \delta_{x,y}
 - d^{\ell}_{1} \delta_{x,y+1} - d^{r}_{1} \delta_{x,y-1} & L_- < x < L_+ \\
d^{\ell}_{1} 
\left(\delta_{y,L_+} - \delta_{y,L_+-1}\right) & x = L_+
\ea \right.
\ee
Comparing with the intensity matrix $Q$ of the biased random walk
of \ref{Rem:RW} one finds $G_{xy}=Q_{xy}$.
Thus with $H = - W^T$ and 
$S_{x\bfeta} := \pi^{-1}(x) R(x,\bfeta) = \mu^{x}(\bfeta)$ one gets
$SW = QS$ which according to \eref{intdualitydef} in the proof of Theorem \ref{Theo:measure2} is equivalent to \eref{TRD1shock} due to reversibility of the
shock exclusion process w.r.t. the measure $\pi(x) = d_1^x$ of the biased random walk.
\qed

\subsubsection{Proof of Theorem \ref{Theo:Evolution1shock}}

The assertion \eref{Evolution1shock} as such, i.e., without specifying the transition
probability, is an immediate corollary of the time-reversed duality 
established in Theorem \eref{Theo:TRD1shock}.
The specific formula \eref{condprobrw} follows from the fact that
the intensity matrix for a biased simple random walk with constant rates
is a tridiagonal Toeplitz matrix. For readers not familiar with such matrices 
we outline how to arrive at \eref{condprobrw}:
(a) The transition probability is by definition the matrix element 
$P(y,t|x,0) = (\exp{G t})_{x,y}$ of the exponential of the intensity matrix.
(b) The non-symmetric intensity matrix can be symmetrized by the 
ground state transformation $\tilde{G} = \hat{\pi} G \hat{\pi}^{-1}$ where $\hat{\pi}$ is the
diagonal matrix with elements $d_1^{x-L_-}$ on the diagonal. (c)
Diagonalizing the tridiagonal Toeplitz matrix $G$ (which is textbook) 
then yields the eigenvalues
$\epsilon_0 =0$ and $\epsilon_p$ \eref{epsrw}
 for $1\leq p \leq L-1$ and the associated eigenvectors of $\tilde{G}$ 
 with components that can be expressed as
\bea 
\Psi_0(y) & = & \sqrt{\frac{d_1^2-1}{d_1^{2L}-1}} d_1^{y-L_-} \\
\Psi_p(y) & = & \sqrt{\frac{2}{L}} 
\frac{\psi_p(y)}
{d_1 - \rme^{-\frac{i \pi p}{L}} } 
\eea
and satisfy the orthogonality relations 
\be 
\sum_{y= -L_-}^{L_+} \Psi_p(y) \Psi^\ast_q(y) = \delta_{p,q} \quad \forall p,q\in \{0,1,\dots,L-1\} .
\ee
Orthogonality can be proved by expanding for $d_1\neq 1$ the denominator in the power series
\be
\frac{1}{d_1 - \rme^{-\frac{i \pi p}{L}}}
= \left\{\ba{ll} \displaystyle
\frac{1}{1-d_1^{-2L}} \sum_{n=0}^{2L-1} d_1^{-n-1} \rme^{-\frac{i n\pi p}{L}} & d_1 > 1 \\
\displaystyle
- \frac{1}{1 - d_1^{2L}} \sum_{n=0}^{2L-1} d_1^{n} \rme^{\frac{i (n+1) \pi p}{L}} & d_1 < 1
\ea\right.
\ee
and using the orthogonality relations
\bea 
\frac{2}{L} \sum_{y=L_-}^{L_+} \sin{\frac{\pi p (y+1-L_-)}{L} }\sin{\frac{\pi q (y+1-L_-)}{L}}  
& = & \sum_{k\in\Z} \left(\delta_{p-q,2kL} - \delta_{p+q,2kL}
\right) \\
\frac{2}{L} \sum_{p=1}^{L-1} \sin{ \frac{\pi p (x+1-L_-)}{L} }\sin{\frac{\pi p (y+1-L_-)}{L}}  
& = & \sum_{k\in\Z} \left(\delta_{x-y,2kL} - \delta_{x+y,2kL}
\right)
\eea
and similar addition formulae involving three trigonometric functions.
(d) It then follows that $P(y,0|x,0) = \delta_{y,x}$
and eigenvector decomposition of the matrix element $(\exp{G t})_{x,y}$
leads to \eref{condprobrw}. Finally, we note that existence of the limit
$t\to\infty$ can be read off the fact that the eigenvalues \eref{epsrw} of the intensity matrix are strictly
positive for all $p \in \{1,\dots,L-1\}$. Uniqueness of the invariant measure
follows from ergodicity of the process.
\qed

\subsection{Reverse duality involving $N$ shocks}

The strategy of the proof is analogous to the case of one shock which
can be used as it is for shock measures that do not have neighbouring
shock positions.

\subsubsection{Proof of Theorem \ref{Theo:TRDNshock}}

\paragraph{Isolated shocks}
If all shocks are isolated, i.e., there is at least one site with density
$\rho_{i}$ between consecutive shocks with densities $\rho^\star_{i}$
$\rho^\star_{i+1}$, then according to (A) and (B) of the projection 
lemma \ref{Lem:hbulk} the generator $H$ acts on the shock measure by acting
on all shocks separately as in the case of a single shock. In other
words, for $N$ shocks at positions $x_i$ such that $x_{i+1} - x_{i} > 1$ for
all $i$ one gets
\bea 
H \ket{\mu^{\bfx}} & = & \left(h^-_{L_-} 
+ \sum_{i=1}^N \left(\tilde{h}_{x_{i-1}} + \tilde{h}_{x_{i}}\right)
+ h^+_{L_+} \right) \ket{\mu^{\bfx}} \\
& = & \sum_{i=1}^N \left[
d^{\ell}_{i} \left(1-\delta_{x_1,L_-} \right) 
\left(\ket{\mu^{\bfx}} - \ket{\mu^{\bfx^{i-}}}\right) \right. \nonumber \\
& & \left. + d^{r}_{i} \left( 1-\delta_{x_N,L_+} \right) \left(\ket{\mu^{\bfx}} - 
\ket{\mu^{\bfx^{i+}}}\right) 
\right]
\eea
which implies invariance of the subspace with $N$ shocks with bulk
densities satisfying the microscopic stability condition \eref{shockcondi}
and arbitrary shock densities provided that
no shocks at nearest neighbour sites appear in the sum over the shock
positions. 

The new situation for which closure needs to be proved
are bunches of neighbouring shocks. We call $n$ shocks
on consecutive sites $x_{k},x_{k+1},\dots,x_{k+n-1}$ with shock
densities $\rho^\star_{i+1},\rho^\star_{i+2},\dots,\rho^\star_{i+n}$
and bulk densities $\rho_{i}$ on site $x_{k-1}$ and $\rho_{i+n}$
on site $x_{k+n}$ satisfying \eref{Nshockcond} with $N=n$ a shock bunch of size $n$ at position $x_k$. An isolated
shock is called a bunch of size 1.
Evidently, for a total of $N$ shocks, the bunch size is in the range
$1\leq n \leq N$, the shock index is in range $0\leq i \leq N-n$, 
and the bunch position is in the range $L_- \leq k \leq L_+ - N+1$.
The new technical ingredient to prove closure for bunches is assertion
(C) in Lemma \ref{Lem:hbulk}.

\paragraph{Bulk bunches of size $n$:} For fixed $n$ we define the product vectors
\bea
\ket{\varXi_{i}} & := & \ket{\rho_{i}\rho^{\star}_{i+1}\rho^{\star}_{i+2}\dots
\rho^{\star}_{i-1+n}\rho^{\star}_{i+n}\rho_{i+n}} \\
\ket{\varOmega_{i}} & := & (\tilde{h}_{12} + \sum_{k=2}^{n} \tilde{h}_{kk+1} + \tilde{h}_{n+1,n+2}) \ket{\varXi_{i}} \\
\ket{\varXi^{\ell}_{i+1}} & := & \ket{\rho^{\star}_{i+1}\rho_{i+1}\rho^{\star}_{i+2}\dots\rho^{\star}_{i-1+n}\rho^{\star}_{i+n}\rho_{i+n}} \\
\ket{\varXi^{r}_{i+n}} & := & \ket{\rho_{i}\rho^{\star}_{i+1}\rho^{\star}_{i+2}\dots \rho^{\star}_{i-1+n}\rho_{i-1+n}\rho^{\star}_{i+n}} \\
\ket{\varUpsilon_{ij}} & := & \ket{\rho_{i}\rho^{\star}_{i+1}\dots\rho^{\star}_{j-1}c_{j}\rho^{\star}_{j+1}\dots\rho^{\star}_{i-1+n}\rho_{i+n}}, \quad i < j \leq i+n \\
\ket{\tilde{\varUpsilon}_{ij}} & := & \ket{\rho_{i}\rho^{\star}_{i+1}\dots\rho^{\star}_{j-1}\tilde{c}_{j}\rho^{\star}_{j+1}\dots\rho^{\star}_{i-1+n}\rho_{i+n}}, \quad i < j \leq i-1+n .
\eea
To compute $\ket{\varOmega_{i}}$ we use Lemma \ref{Lem:hbulk}, specifically assertion (A)
for the action of $\tilde{h}_{12}$, assertion (B)
for the action of $\tilde{h}_{n+1,n+2}$, and assertion (A)
for the action of the remaining local bulk jump operators $\tilde{h}_{kk+1}$ in the summation symbol. 
Setting $a^{(1)}=\rho_{i}$, $\tilde{a}^{(1)}=\rho_{i+1}$, $b^{(1)}=\rho^{\star}_{i+1}$,
$\tilde{b}^{(1)}=c_{i+1}$ in (A),
$a = \rho^{\star}_{k}$, $\tilde{a} = \rho^{\star}_{k+1}$, $\tilde{c}=\tilde{c}_k$, $c=c_{k+1}$ 
for $i+1 \leq k \leq i+n-1$ in (C), and $a^{(2)}=\rho_{i+n-1}$, $\tilde{a}^{(2)} = \rho_{i+n}$, $b^{(2)} = \rho^{\star}_{i+n}$,
$\tilde{b}^{(2)} = \tilde{c}_{i+n}$
this yields
\bea 
\ket{\varOmega_{i}} & = & 
\left(r+\ell 
+ \sum_{k=i+1}^{i+n} \left(\tilde{d}(\rho^{\star}_{k},\tilde{c}_{k}) 
- \tilde{d}(\rho^{\star}_{k},c_{k})\right)
\right) \ket{\varXi_{i}} \nonumber \\
& & - d(\rho_{i},\rho_{i+1}) \ket{\varXi^{\ell}_{i+1}} + d(\rho_{i+n}, \rho_{i+n-1}) \ket{\varXi^{r}_{i+n}}
 \nonumber \\
& & + \sum_{k=i+1}^{i+n}
\left(d(\rho^{\star}_{k},c_{k}) \ket{\varUpsilon_{ik}}  - 
d(\rho^{\star}_{k},\tilde{c}_k) \ket{\tilde{\varUpsilon}_{ik}}\right) 
\eea
Without loss of generality we choose $c_{k} = \tilde{c}_{k} \neq \rho^{\star}_{k}$ so that
the sums vanish and using \eref{vi} and \eref{wi} one finds
\bea 
\ket{\varOmega_{i}} 
& = & \left(r+\ell \right) \ket{\varXi_{i}} \nonumber \\
& & - d^{\ell}_{i+1} \ket{\varXi^{\ell}_{i+1}} - d^{r}_{i+n} \ket{\varXi^{r}_{i+n}}
\nonumber \\
& = & \left[d^{\ell}_{i+1} + d^{r}_{i+1} + (r-\ell)(\rho_{i} - \rho_{i-1})\right] \ket{\varXi_{i}} \nonumber \\
& & - d(\rho_{i},\rho_{i+1}) \ket{\varXi^{\ell}_{i+1}} + d(\rho_{i+n}, \rho_{i+n-1}) \ket{\varXi^{r}_{i+n}} 
\eea
if and only if \eref{shockcondi} which is equivalent to \eref{hab2}
holds for all consecutive shock densities and all consecutive bulk 
densities.

Applying this relation to  a measure $\ket{\mu^{\bfx}}$ with $N$ shocks 
at bulk positions $L_- < x_1 <\dots <x_N < L_+$ arranged in $k$ 
bunches one thus finds an invariant subspace under the action of the 
generator by noting
\bea 
H \ket{\mu^{\bfx}} & = & [\epsilon_-(\rho_0) + \epsilon_+(\rho_N)] \ket{\mu^{\bfx}} \nonumber \\
& & + \sum_{i=1}^{N} (1-\delta_{x_{i-1},x_{i}-1})
\left( [d^{\ell}_{i} - (r-\ell)\rho_{i-1}] \ket{\mu^{\bfx}}
- d^{\ell}_{i} \ket{\mu^{\bfx^{i-}}} \right) \nonumber \\
& & + \sum_{i=1}^{N} (1-\delta_{x_{i+1},x_{i}+1})
\left([d^{r}_{i} + (r-\ell)\rho_{i}] \ket{\mu^{\bfx}}
- d^{r}_{i} \ket{\mu^{\bfx^{i+}}}\right) \nonumber \\
& = & \sum_{i=1}^{N} d^{\ell}_{i} (1-\delta_{x_{i-1},x_{i}-1}) 
\left(\ket{\mu^{\bfx}} - \ket{\mu^{\bfx^{i-}}}\right) \nonumber \\
& & + \sum_{i=1}^{N} d^{r}_{i} (1-\delta_{x_{i+1},x_{i}+1})
\left(\ket{\mu^{\bfx}} - \ket{\mu^{\bfx^{i+}}}\right) \nonumber \\
& & + (r-\ell) 
\left(\delta_{x_{1},L_-} \rho_{0} - \delta_{L_+,x_{N}} \rho_{N} \right)
\ket{\mu^{\bfx}} 
\eea
where in the first equality the bunch property and in the second equality 
we have used cancellation of the terms in the telescopic sum over the 
diagonal density terms with the boundary eigenvalues. Observing finally 
that $\delta_{x_{1},L_-} = \delta_{L_+,x_{N}} = 0$ for $L_- < x_1 
<\dots <x_N < L_+$ we arrive at
\bea 
H \ket{\mu^{\bfx}} 
& = & \sum_{i=1}^{N} d^{\ell}_{i} (1-\delta_{x_{i-1},x_{i}-1}) 
\left(\ket{\mu^{\bfx}} - \ket{\mu^{\bfx^{i-}}}\right) \nonumber \\
& & + \sum_{i=1}^{N} d^{r}_{i} (1-\delta_{x_{i+1},x_{i}+1})
\left(\ket{\mu^{\bfx}} - \ket{\mu^{\bfx^{i+}}}\right) \\
& = & \sum_{i=1}^{N} \left[w^{\ell}_{i} 
\left(\ket{\mu^{\bfx}} - \ket{\mu^{\bfx^{i-}}}\right) + w^{r}_{i} \left(\ket{\mu^{\bfx}} - \ket{\mu^{\bfx^{i+}}}\right) \right]
\label{Hmx}
\eea
for shocks not at the boundary sites. It remains to prove closure
when one or two shocks are at the boundary sites $L_\pm$.

\paragraph{Bunch of shocks at the boundaries:}
We first consider a bunch of $n$ shocks at the left boundary but no shock at the right boundary, i.e., $x_1 = L_-$ and $x_N < L_+$. For fixed $n$ we define the product vectors
\bea
\ket{\varXi^{-}} & := & \ket{\rho^{\star}_{1}\rho^{\star}_{2}\dots
\rho^{\star}_{n-1}\rho^{\star}_{n}\rho_{n}} \\
\ket{\varOmega^{-}} & := & (h^-_1 + \sum_{k=1}^{n-1} \tilde{h}_{kk+1} + \tilde{h}_{n,n+1}) \ket{\varXi_{-}} \\
\ket{\varXi^{-}_{n}} & := & \ket{\rho^{\star}_{1}\rho^{\star}_{2}\dots \rho^{\star}_{n-1}\rho_{n-1}\rho^{\star}_{n}} \\
\ket{\varUpsilon^{-}_{j}} & := & \ket{\rho^{\star}_{1}\dots\rho^{\star}_{j-1}c_{j}\rho^{\star}_{j+1}\dots\rho^{\star}_{n-1}\rho^{\star}_{n}\rho_{n}}, \quad 1 \leq j \leq n \\
\ket{\tilde{\varUpsilon}^{-}_{j}} & := & \ket{\rho^{\star}_{1}\dots\rho^{\star}_{j-1}\tilde{c}_{j}\rho^{\star}_{j+1}\dots\rho^{\star}_{n-1}\rho^{\star}_{n}\rho_{n}}, \quad 1 \leq j \leq n \\
\eea
To compute $\ket{\varOmega^{-}}$ we use assertion (B) in
Lemma \ref{Lem:hbulk}
for the action of $\tilde{h}_{n,n+1}$ and assertion (C)
for the action of the remaining local bulk jump operators $\tilde{h}_{kk+1}$ in the summation symbol. With $a=\rho^\star_{1}$ and $b=c_1$
for the boundary terms, 
$a=\rho^\star_{k}$, $\tilde{a} = \rho^\star_{k+1}$, $c=c_{k+1}$,
and $\tilde{c}=\tilde{c}_k$ for assertion (C), and $b^{(2)} = 
\rho^\star_{n}$, $\tilde{a}^{(2)} = \rho_{n}$, $a^{(2)} = \rho_{n-1}$,
$\tilde{b}^{(2)} = \tilde{c}_n$ for assertion (B)
this yields
\bea 
\ket{\varOmega^{-}} 
& = &
\left(\tilde{d}_-(\rho^\star_{1},c_1) - \tilde{d}(\rho^\star_{1},c_1)\right) \ket{\varOmega^{-}}  \nonumber \\
& & - \left(d_-(\rho^\star_{1},c_1) - d(\rho^\star_{1},c_1)\right)  \ket{\varUpsilon^{-}_{1}} \nonumber \\
& & + \sum_{k=1}^{n-1}
\left(\tilde{d}(\rho^\star_{k},\tilde{c}_k) - \tilde{d}(\rho^\star_{k+1},c_{k+1})\right) \ket{\varOmega^{-}}
\nonumber \\
& & - \sum_{k=1}^{n-1} d(\rho^\star_{k},\tilde{c}_k) \ket{\tilde{\varUpsilon}^{-}_{k}}
+ \sum_{k=1}^{n-1} d(\rho^\star_{k+1},c_{k+1}) \ket{\varUpsilon^{-}_{k+1}} \nonumber \\
& & + \left(r + \tilde{d}(\rho^\star_{n}, \tilde{c}_n) \right) \ket{\varOmega^{-}} \nonumber \\
& & + d(\rho_{n}, \rho_{n-1}) \ket{\varXi^{-}_{n}}
- d(\rho^\star_{n}, \tilde{c}_n) \ket{\tilde{\varUpsilon}^{-}_{n}}
\nonumber \\
& = &
\left(r + \tilde{d}_-(\rho^\star_{1},c_1)\right) \ket{\varOmega^{-}}  - d_-(\rho^\star_{1},c_1) \ket{\varUpsilon^{-}_{1}} 
+ d(\rho_{n}, \rho_{n-1}) \ket{\varXi^{-}_{n}}
\nonumber \\
& & 
+ \sum_{k=1}^{n} \left(
d(\rho^\star_{k},c_{k}) \ket{\varUpsilon^{-}_{k}} 
- d(\rho^\star_{k},\tilde{c}_k) \ket{\tilde{\varUpsilon}^{-}_{k}}\right). 
\eea
Now we choose use without loss of generality 
$\tilde{c}_k = c_k$ and observe
that by assumption $\rho^\star_1=\alpha/(\alpha+\gamma)$ so that
$d_-(\rho^\star_{1},c_1) = \tilde{d}_-(\rho^\star_{1},c_1) = 0$ 
for any choice of $c_1$. Therefore 
\bea 
\ket{\varOmega^{-}} 
& = &
r \ket{\varOmega^{-}}  
+ d(\rho_{n}, \rho_{n-1}) \ket{\varXi^{-}_{n}} \nonumber \\
& = &
[d^{r}_n + (r - \ell) \rho_n] \ket{\varOmega^{-}}
- d^{r}_n \ket{\varXi^{-}_{n}}
\eea
where in the second equality \eref{vi} was used and one gets
for consecutive bunches with a total of $N$ shocks
\bea 
H \ket{\mu^{\bfx}} & = & \sum_{i=1}^{N} d^{\ell}_{i} (1-\delta_{x_{i-1},x_{i}-1}) 
\left(\ket{\mu^{\bfx}} - \ket{\mu^{\bfx^{i-}}}\right) \nonumber \\
& & + \sum_{i=1}^{N} d^{r}_{i} (1-\delta_{x_{i+1},x_{i}+1})
\left(\ket{\mu^{\bfx}} - \ket{\mu^{\bfx^{i+}}}\right) \nonumber \\
& & - (r-\ell) \rho_{0}
\left(1 - \delta_{x_{1},L_-}\right)
\ket{\mu^{\bfx}} .
\eea
Since for a shock at the left boundary site one has $1 - \delta_{x_{1},L_-} = 0$ one recovers \eref{Hmx}
also for $x_1 = L_-$. For a bunch of shocks only at the right boundary
or two bunches on both boundaries similar computations using \eref{bNshockcond2} with $N=m+1$ also lead to
\eref{Hmx} which is therefore valid for any set of sites $\bfx$.

Hence invariance of the subspace spanned by the vectors $\ket{\mu^{\bfx}}$ is proved. The closure relation \eref{Hmx}
can be recast as
\be 
H \ket{\mu^{\bfx}} = - \sum_{\bfy} Q_{\bfx\bfy} \ket{\mu^{\bfy}}
\label{Hmx2}
\ee
with the matrix elements $Q_{\bfx\bfy}$ \eref{shockasepQ}
of the intensity matrix for $N$-particle shock exclusion process. With the same arguments as at the end of the proof of Theorem \ref{TRD1shock} one arrives at \eref{TRDNshock}.
\qed

\subsubsection{Proof of Theorem \ref{Theo:Nshockevolution}}

The time evolution formula \eref{shockevolutionN}
follows from Theorem \ref{Theo:measure2}.
Existence of the limit \eref{invmeasN}
follows trivially from the finite state space of the ASEP defined on $\Lambda_L$ (which guarantees existence of an invariant measure)
and uniqueness follows from ergodicity which is guaranteed
by the strict positivity of the boundary rates.
\qed

\section*{Acknowledgements}

GMS thanks the Galileo-Galilei Institute in Florence for kind hospitality
and support during the scientific program on ``Randomness, Integrability, and 
Universality'' and particularly C. Giardin\`a and J. de Gier for inspriring
discussions and pointing out useful references. This work is financially 
supported by FCT/Portugal through CAMGSD, IST-ID, Projects UIDB/04459/2020 
and UIDP/04459/2020. 

\appendix
\renewcommand{\theequation}{\Alph{section}.\arabic{equation}}

\section{The XXZ quantum spin chain with non-diagonal boundary fields}
\label{App:A}
\setcounter{equation}{0}

We denote by $\mathds{1}$ the two-dimensional unit matrix and by 
$\mathbf{1}$ the unit matrix of dimension $2^L$. We also recall the 
standard definitions the Pauli matrices and some of their linear combinations.
\begin{df}
\label{Def:localmatrices}
The Pauli matrices are the $2\times2$ matrices
\be 
\sigma^x = \left( \ba{cc} 0 & 1 \\ 1 & 0 \ea \right), \quad
\sigma^y = \left( \ba{cc} 0 & -i \\ i & 0 \ea \right), \quad
\sigma^z = \left( \ba{cc} 1 & 0 \\ 0 & -1 \ea \right)
\ee
and
\be 
\sigma^\pm = \half \left( \sigma^x \pm i \sigma^y \right)
\ee
are the particle annihilation and creation operators. The diagonal
matrices
\be
\hat{n} = \half \left( \mathds{1} - \sigma^z \right), \quad 
\hat{v} = \half \left( \mathds{1} + \sigma^z \right)
\ee
are the particle and vacancy projectors respectively.
\end{df}

With these notations the negative transpose $H$ \eref{generator}
of the intensity matrix $M$ \eref{MASEP}  is given by
\bea
H & = & - \sum_{k=L_-}^{L_+-1}[ r (\sigma^+_{k} \sigma^-_{k+1}
- \hat{n}_{k} \hat{v}_{k+1}) + \ell (\sigma^-_{k} \sigma^+_{k+1}
- \hat{v}_{k} \hat{n}_{k+1})] \nonumber \\
& & - [\alpha(\sigma^-_{L_-} - \hat{v}_{L_-}) 
+ \gamma (\sigma^+_{L_-} - \hat{n}_{L_-}) + 
\delta(\sigma^-_{L_+} - \hat{v}_{L_+}) 
+ \beta 
(\sigma^+_{L_+} - \hat{n}_{L_+})] 
\eea
A transformation with the diagonal matrix
\be 
\hat{Q} := \prod_{k=L_-}^{L_+-1} q^{k \hat{n}_k}
\ee
symmetrizes the bulk jump matrices. With the parameters
\bea 
& & \theta := \ln{q}, \quad w := \sqrt{r\ell} 
\label{thetadef} \\
& & E_0 := (L-1) \cosh{(\theta)} + \frac{\alpha+\beta+\gamma+\delta}{w}
\label{E0def}
\eea
and with the parametrization
\bea 
\alpha & = & \frac{w}{2} \sinh{\theta}
\frac{\rme^{\phi_{-} - \psi_-}}{\sinh{\phi_{-}} \cosh{\psi_{-}}} 
\label{apar} \\
\gamma & = & \frac{w}{2} \sinh{\theta} 
\frac{\rme^{\psi_{-} - \phi_{-} }}{\sinh{\phi_{-}} \cosh{\psi_{-}}} 
\label{gpar} \\
\beta & = & \frac{w}{2} \sinh{\theta}
\frac{\rme^{\phi_{+} - \psi_{+}}}{\sinh{\phi_{+}} \cosh{\psi_{+}}}  \label{bpar} \\
\delta & = & \frac{w}{2} \sinh{\theta}
\frac{\rme^{- \phi_{+} + \psi_{+}}}{\sinh{\phi_{+}} \cosh{\psi_{+}}} 
\label{dpar} \\
\theta_- & = & \psi_{-} - \phi_{-} + \theta L_- 
\label{thetamdef} \\
\theta_+ & = & \phi_{+} - \psi_{+} + \theta L_+
\label{thetapdef}
\eea

one gets the quantum Hamiltonian 
\bea
H^{XXZ}  & = & \hat{Q}^{-1} H \hat{Q} \\
& = & - \frac{w}{2} \left\{\sum_{k=L_-}^{L_+-1} \left[ \sigma^x_{k} \sigma^x_{k+1} + \sigma^y_{k} \sigma^y_{k+1}
+ \cosh{(\theta)}
\sigma^z_{k} \sigma^z_{k+1}\right] 
- E_0 \mathbf{1} \right. \nonumber \\
& & + \frac{\sinh{\theta}}{\sinh{\phi_{-}} \cosh{\psi_{-}}}  \left(\rme^{-\theta_-} \sigma^-_{L_-} 
+ \rme^{\theta_-} \sigma^+_{L_-} 
+ \sinh{\psi_{-}} \cosh{\phi_{-}} \sigma^z_{L_-} \right)
\nonumber \\
& & \left. + 
\frac{\sinh{\theta}}{\sinh{\phi_{+}} \cosh{\psi_{+}}} 
\left(\rme^{- \theta_{+}} \sigma^-_{L_+} 
+ \rme^{\theta_+} \sigma^+_{L_+}
- \sinh{\psi_{+}} \cosh{\phi_{+}} \sigma^z_{L_+} \right)
\right\} .
\label{HXXZ}
\eea
of the XXZ spin chain with non-diagonal
and non-hermitian boundary fields in the form presented in \cite{Nepo03b} for $w=1$ and $E_0=0$.


In terms of the parameters \eref{thetadef} - \eref{thetapdef}
the boundary functions \eref{kappadef}
are given by
\bea 
\kappa_+(\alpha,\gamma) & = & \rme^{\phi_{-} +\psi_{-}-\theta L_- + \theta_-}\\
\kappa_+(\beta,\delta) & = & \rme^{\phi_+ + \psi_+ - \theta_+ + \theta L_+}
\eea
so that the condition \eref{bNshockcond} reads
\be 
\phi_{-} +\psi_{-} + \phi_+ + \psi_+ = \theta_+ - \theta_- + (2N - L +1) \theta .
\ee
This is the integrability condition found in \cite{Nepo03b}.
\qed
The special condition \eref{bNshockcond2} reads
\bea
\frac{\alpha\beta}{\gamma\delta} = \rme^{2(\phi_{-} + \phi_{+} - \psi_{-} - \psi_{-})} = \rme^{-2 \theta (N-M)}
\eea
which can be recast as
\bea
\theta_{+} - \theta_{-} + \theta (N-M-L+1) = 0.
\eea

\section{On the parametrization of the open ASEP}
\label{App:B}
\setcounter{equation}{0}

With the parametrization \eref{lbr} and \eref{rbr} one gets
\bea 
& & \kappa_+(\alpha,\gamma) = z_-^{-1}, \quad 
\kappa_-(\alpha,\gamma) = - \frac{\ell+\omega_-}{r+\omega_-} 
\label{kag} \\
& & \kappa_+(\beta,\delta) = z_+, \quad 
\kappa_-(\beta,\delta) = - \frac{\ell+\omega_+}{r+\omega_+}.
\label{kbd}
\eea
The independence of the function $\kappa_+(\alpha,\gamma)$
of $\omega_-$ can be seen by noting that 
$x/\rho - y/(1-\rho) = r-\ell$ and therefore
$y - x + r - \ell = x z_-^{-1} - y z_-$ so that the argument of the square
root becomes $(x z_-^{-1} + y z_-)^2$. With the role of $\rho$ and $1-\rho$ interchanged one then also realizes that $\kappa_+(\beta,\delta)$
is independent of  $\omega_+$. Thus on the parameter manifold given by
\be
\label{aga}
(r - \ell)\rho_-(1-\rho_-) = \alpha(1-\rho_-) - \gamma \rho_-
\ee
one has $\kappa_+(\alpha,\gamma) = z_-^{-1}$
and on the manifold
\be
\label{bdb}
(r - \ell)\rho_+(1-\rho_+) = \beta \rho_+ - \delta(1-\rho_+)
\ee
one has $\kappa_+(\beta,\delta) = z_+$.

For $q>1$, boundary densities $\rho_\pm \in (0,1)$, and boundary 
parameters in $\R^+$ one has $\omega_\pm > - \ell$
and correspondingly $-1 < \kappa_-(x,y) < 0$. Moreover, 
$\kappa_+(x,y)\in\R^+$. Then \eref{kpm} implies $\rho_0<\sigma^\star_1<\rho_1$.
On the other hand, for $q<1$ where a downward shock with $\rho_0 >\rho_1$ is macroscopically stable one gets $-\infty < \kappa_-(x,y) < -1$ and 
therefore $\rho_1<\sigma^\star_1<\rho_0$. Hence either way the shock
density is between the boundary densities of the shock.

Alternatively, \eref{bNshockcond2} can be written
\be 
\frac{\alpha\beta}{\gamma\delta} = q^{2(N-M)} 
\label{bNshockcond5}
\ee
In terms of the parameters $\rho_\pm$ the manifold
$\mathcal{B}_N$ defined by \eref{bNshockcond} is given by
\be 
\frac{z_+}{z_-} = q^{2N}, \quad N \in \N^+ 
\label{bNshockcond3}
\ee
and \eref{bNshockcond2} defining the submanifold $\mathcal{B}^M_N$
can be written as the constraint
\be
\frac{(r+\omega_-) (r+\omega_+)}{(\ell+\omega_-)(\ell+\omega_+)} = q^{2M}, \quad
1 \leq M \leq N
\label{bNshockcond4}
\ee
on the boundary parameters $\omega_\pm$. 

%
%
%


\begin{thebibliography}{99}%

\bibitem{Alca94}
F.C. Alcaraz, M. Droz, M. Henkel, and V. Rittenberg,
Reaction-Diffusion Processes, Critical Dynamics and Quantum ak16ains,
Ann. Phys. \textbf{230}, 250--302 (1994).

\bibitem{Baha12}
C. Bahadoran,
Hydrodynamics and Hydrostatics for a Class of Asymmetric Particle Systems with Open Boundaries,
Commun. Math. Phys. \textbf{310}(1), 1--24 (2012).

\bibitem{Bala04}
Bal\'azs, M.: Multiple shocks in bricklayers’ model. J. Stat. Phys. 117(1–2), 77–98 (2004). 

\bibitem{Bala10}
Bal\'azs, M., Farkas, G., Kov\'acs, P., and R\'akos, A.: 
Random walk of second class particles in product shock measures.
J. Stat. Phys. \textbf{139}(2), 252--279 (2010)

\bibitem{Bala19}
Bal\'azs, M., Duffy, L. and Pantelli, D. q-Zero Range has Random Walking Shocks. J Stat Phys 174, 958--971 (2019). 

%
\bibitem{Beli02} 
Belitsky, V., Sch\"utz, G.M.: 
Diffusion and scattering of shocks in the partially asymmetric simple exclusion process.
Electron. J. Probab. \textbf{7}, paper 11, 1-21 (2002)

%
\bibitem{Beli13} 
V.Belitsky and G.M. Sch\"utz,
Microscopic structure of shocks and antishocks in the ASEP conditioned on low current,
J. Stat. Phys. \textbf{152}, 93--111 (2013).

\bibitem{Beli13b} 
V Belitsky and G M Sch\"utz, 
Antishocks in the ASEP with open boundaries conditioned on low current
J. Phys. A: Math. Theor. \textbf{46}, 295004 (2013)

%
%
\bibitem{Beli15b}
V. Belitsky and G.M. Sch\"utz,  
Self-Duality for the Two-Component Asymmetric Simple Exclusion Process.
J. Math. Phys. \textbf{56}, 083302 (2015).

\bibitem{Beli18}
V. Belitsky and G. M. Sch\"utz,
Self-duality and shock dynamics in the $n$-species priority ASEP,
Stoch. Proc. Appl. \textbf{128}, 1165--1207 (2018).

%
%
%
%
%

\bibitem{Boro20}
A. Borodin, I. Corwin,
Dynamic ASEP, Duality, and Continuous $q^{-1}$-Hermite Polynomials,
Int. Math. Res. Notices, Vol. 2020, No. 3, 641--668 (2020).

\bibitem{Bryc19}
W. Bryc and M. Swieca,
On Matrix Product Ansatz for Asymmetric Simple Exclusion Process with Open Boundary in the Singular Case,
J. Stat. Phys.  \textbf{177}, 252--284 (2019)

\bibitem{Cari13}
G. Carinci, C. Giardina, C. Giberti, and F. Redig,
Duality for Stochastic Models of Transport,
J. Stat. Phys.  \textbf{152}, 657--697 (2013).

%
%

%
%
%
\bibitem{Cari16}
G. Carinci, C. Giardin\`a, F. Redig, and T. Sasamoto, 
A generalized Asymmetric Exclusion Process with $U_q(\mathfrak{sl}_2)$ 
stochastic duality, 
Probab. Theory Relat. Fields \textbf{166}, 887--933 (2016).

\bibitem{Cari19}
G. Carinci, C. Franceschini, C. Giardin\`a, W. Groenevelt, and F. Redig,
Orthogonal Dualities of Markov Processes and Unitary Symmetries,
SIGMA \textbf{15}, 053 (2019)

\bibitem{Cari21}
G. Carinci, C. Franceschini, and W. Groenevelt,
$q$-Orthogonal dualities for asymmetric particle systems
Electron. J. Probab. \textbf{26}: 1--38 (2021). 

\bibitem{Chak16}
S. Chakraborty, S. Pal, S. Chatterjee and M. Barma,
Large compact clusters and fast dynamics in coupled nonequilibrium systems,
Phys. Rev. E \textbf{93}, 050102(R) (2016).

\bibitem{Chen20}
Z. Chen, J. de Gier, and M. Wheeler,
Integrable Stochastic Dualities and the Deformed Knizhnik–Zamolodchikov Equation,
Int. Math. Res. Notices, Vol. 2020, No. 19, 5872--5925 (2020).

\bibitem{Chun16}
Chung, K. L., Walsh, J. B., Markov processes, 
Brownian motion, and time symmetry. 2nd edition, 
(Springer, New York, 2005)

\bibitem{Clin03} 
M. Clincy and M.R. Evans,
Phase transition in the ABC model,
Phys. Rev. E \textbf{67}, 066115 (2003).

%

%

\bibitem{Cram10}
N. Cramp\'e, E. Ragoucy, and D. Simon,
Construction of a coordinate Bethe ansatz for the asymmetric simple exclusion process with open boundaries,
J. Stat. Mech. \textbf{2010}, P11038 (2010)

\bibitem{deGi05}
J. de Gier and F.H.L Essler, 
Bethe ansatz solution of the asymmetric exclusion process with open boundaries, 
Phys. Rev. Lett. \textbf{95}, 240601 (2005).

\bibitem{deGi06}
J. de Gier and F. H. L. Essler,
Exact spectral gaps of the asymmetric exclusion process with open boundaries,
J. Stat. Mech. \textbf{2006}, P12011 (2006).

\bibitem{Derr93}
B. Derrida, S.A. Janowsky, J.L. Lebowitz, and E.R. Speer,
Exact solution of the totally asymmetric simple exclusion process: Shock profiles, 
J. Stat. Phys. \textbf{73}, 813--842 (1993).

\bibitem{Derr93a}
B.~Derrida, M.R.~Evans, V.~Hakim, and V.~Pasquier,
Exact solution of a 1D asymmetric exclusion model using a matrix formulation,
J. Phys. A: Math. Gen. \textbf{26}, 1493--1517 (1993).

\bibitem{Derr97}
B. Derrida, J. L. Lebowitz, E. R. Speer, 
Shock profiles in the asymmetric simple exclusion process in one dimension, 
J. Stat. Phys. \textbf{89}, 135--167 (1997).

\bibitem{deVe94}
H.-J. de Vega and A. Gonzalez-Ruiz, 
Boundary K-matrices for the XYZ, XXZ and XXX spin chains,
J. Phys. A: Math. Gen. \textbf{27}, 6129--6137 (1994).

\bibitem{Dudz00}
M. Dudzi\'nski and G.M. Sch\"utz,
Relaxation spectrum of the asymmetric exclusion process with open boundaries.
J. Phys. A: Math. Gen. \textbf{33}, 8351--8364 (2000).

\bibitem{Essl96}
F.H.L. Essler and V. Rittenberg, 
Representations of the quadratic algebra and partially asymmetric diffusion with open boundaries
J. Phys. A: Math. Gen. \textbf{29}, 3375--3407 (1996). 

\bibitem{Evan98b} 
M.R. Evans,  Y. Kafri, H.M. Koduvely, and D. Mukamel,
Phase separation and coarsening in one-dimensional driven diffusive 
systems: Local dynamics leading to long-range Hamiltonians,
Phys. Rev. E {\bf 58} 2764--2778 (1998).

\bibitem{Ferr91}
Ferrari, P.A., Kipnis, C., Saada, E.: 
Microscopic Structure of Travelling Waves in the Asymmetric Simple Exclusion Process. 
Ann. Probab., \textbf{19}(1) 226--244 (1991)

\bibitem{Ferr94} 
P.A. Ferrari and L.R.G. Fontes, 
Shock fluctuations in the asymmetric simple exclusion process,
Probab. Theory Relat. Fields \textbf{99}, 305--319 (1994)


\bibitem{Ferr00}
P.A. Ferrari, L.R.G. Fontes, and M.E. Vares, 
The asymmetric simple exclusion model with multiple shocks,
Ann. Inst. H. Poincar\'e Probab. Stat. \textbf{36}(2), 109--126 (2000)


\bibitem{Fras20}
R. Frassek, C. Giardin\`a, and  J. Kurchan,
Duality and hidden equilibrium in transport models,
SciPost Phys. \textbf{9}, 054 (2020)

\bibitem{Giar09}
C. Giardin\`a, J. Kurchan, F. Redig,  and K. Vafayi,
Duality and Hidden Symmetries in Interacting Particle Systems,
J. Stat. Phys. \textbf{135}, 25--55 (2009).

%
%

\bibitem{Henk95}
M. Henkel, E. Orlandini, and G.M. Sch\"utz, 
Equivalences between stochastic processes, 
J. Phys. A: Math. Gen. \textbf{28}, 6335--6344 (1995).

\bibitem{Inam94}
T. Inami and H. Konno, 
Integrable XYZ Spin Chain with Boundaries,
J. Phys. A: Math Gen. {\bf 27}, L913--L918 (1994).

\bibitem{Jafa07}
F. H. Jafarpour and S. R. Masharian,
Matrix Product Steady States as Superposition of Product Shock Measures in 1D Driven Systems
J. Stat. Mech. P10013 (2007)

\bibitem{Jafa09}
F. H. Jafarpour and S. R. Masharian,
Temporal evolution of product shock measures in TASEP with sublattice-parallel update,
Phys. Rev. E \textbf{79}, 051124 (2009) 

\bibitem{Jans14}
S. Jansen and N. Kurt,
On the notion(s) of duality for Markov processes,
Prob. Surveys \textbf{11}, 59--120 (2014).

%
\bibitem{Kafr03}
Y. Kafri, E. Levine, D. Mukamel, G.M. Sch\"utz, and R.D. Willmann
Phase-separation transition in one-dimensional driven models,
Phys. Rev. E \textbf{68}, 035101(R) (2003)

\bibitem{Kolo98b}
A.B. Kolomeisky, G.M. Sch\"utz, E.B. Kolomeisky and J.P. Straley, 
Phase diagram of one-dimensional driven lattice gases with open boundaries,
J. Phys. A: Math. Gen. \textbf{31}, 6911--6919 (1998).

\bibitem{Kreb03}
K. Krebs, F.H. Jafarpour, and G.M. Sch\"utz,
Microscopic structure of travelling wave solutions in a class of stochastic interacting particle systems,
N. J. Phys. \textbf{5}, 145.1--145.14 (2003).

\bibitem{Krug91}
J Krug, 
Boundary-induced phase transitions in driven diffusive systems. 
Phys. Rev. Lett. \textbf{67}, 1882--1885 (1991).

%
\bibitem{Kuan18}
J Kuan,
A Multi-species ASEP(q,j) and q-TAZRP with Stochastic Duality,
Int. Math. Res. Notices, Vol. 2018, No. 17, 5378--5416.

\bibitem{Kuan21}
J Kuan,
Algebraic Symmetry and Self–Duality of an Open ASEP,
Math. Phys. Anal. Geom. \textbf{24} 12 (2021)

\bibitem{Kuan22}
J Kuan, Two Dualities: Markov and Schur–Weyl, 
Int. Math. Res. Notices, Vol. 2022, No. 13, 9633--9662 (2022).

\bibitem{Lahi00}
R. Lahiri, M. Barma, and S. Ramaswamy,
Strong phase separation in a model of sedimenting lattices.
Phys. Rev. E {\bf 61}, 1648--1658 (2000).

%
\bibitem{Ligg76}
T.M. Liggett, 
Coupling the simple exclusion process. 
Ann. Probab. \textbf{4}, 339--356 (1976).

\bibitem{Ligg85}
T.M. Liggett,  
{\it Interacting particle systems}
Springer, Berlin, (1985).

\bibitem{Ligg99} 
T. M.~Liggett, 
{\it Stochastic Interacting Systems: Contact, Voter and Exclusion Processes}
Springer, Berlin (1999).

\bibitem{Lin19}
Y. Lin,
Markov duality for stochastic six vertex model
Electron. Commun. Probab. \textbf{24}, 1-17 (2019). 

\bibitem{Lloy96}
P. Lloyd, A. Sudbury, and P. Donnelly,
Quantum operators in classical probability theory: 
I. ``Quantum spin'' techniques and the exclusion model of diffusion,
Stoch. Proc. Appl. \textbf{61}, 205--221 (1996).

\bibitem{Mall97}
K. Mallick and S. Sandow,
Finite-dimensional representations of the quadratic algebra: Applications to the exclusion process
J. Phys. A: Math. Gen. \textbf{30}, 4513--4526 (1997).

\bibitem{Nepo03b}
R. I. Nepomechie and F. Ravanini,
Completeness of the Bethe Ansatz solution of the open XXZ chain with nondiagonal boundary terms,
J. Phys. A: Math. Gen. \textbf{37}, 433--440 (2004).

\bibitem{Nepo04}
Rafael I. Nepomechie,
Bethe Ansatz solution of the open XXZ chain with nondiagonal boundary terms,
J. Phys. A: Math. Gen. \textbf{37}, 433--440 (2004).

%
\bibitem{Ohku17}
J Ohkubo,
On dualities for SSEP and ASEP with open
boundary conditions, J. Phys. A: Math. Theor. \textbf{50} 095004 (2017).

%

\bibitem{Redi18}
F. Redig, F. Sau,
Stochastic duality and eigenfunctions,
https://arxiv.org/abs/1805.01318

\bibitem{Reza95}
F. Rezakhanlou,
Microscopic structure of shocks in one conservation laws
Annales de l’I. H. P., section C, tome 12, no 2 (1995), p. 119--153

\bibitem{Sand94a}
S. Sandow and G.M. Sch\"utz,
On $U_q[SU(2)]$-Symmetric Driven Diffusion,
Europhys. Lett. \textbf{27}, 7--12 (1994).

%
\bibitem{Sand94b}
S. Sandow, 
Partially asymmetric exclusion process with open boundaries,
Phys. Rev. E \textbf{50}, 2660--2667 (1994). 

\bibitem{Sant02}
Santen, L., Appert, C. The Asymmetric Exclusion Process Revisited: Fluctuations and Dynamics in the Domain Wall Picture. J. Stat. Phys. \textbf{106}, 187--199 (2002). 

\bibitem{Schu93b} 
G. Sch\"{u}tz and E. Domany, 
Phase transitions in an exactly soluble one- dimensional asymmetric exclusion model. 
J. Stat. Phys. \textbf{72}, 277--296 (1993).


\bibitem{Schu94}
G. Sch\"utz, and S. Sandow,
Non-abelian symmetries of stochastic processes: derivation of correlation functions 
for random vertex models and disordered interacting many-particle systems,
Phys. Rev. E \textbf{49}, 2726--2744 (1994).

\bibitem{Schu95}
G.M. Sch\"utz, 
Diffusion-annihilation in the presence of a driving field.
J. Phys. A: Math. Gen. \textbf{28},  3405--3415 (1995).

\bibitem{Schu97} 
G.M. Sch\"utz, 
Duality relations for asymmetric exclusion processes,
J. Stat. Phys. \textbf{86}, 1265--1287 (1997).

\bibitem{Schu97b} 
G.M. Sch\"utz,
The Heisenberg chain as a dynamical model
for protein synthesis - Some theoretical and experimental results,
Int. J. Mod. Phys. B \textbf{11}, 197--202 (1997).


%


\bibitem{Schu01} 
G.M. Sch\"utz, 
Exactly solvable models for many-body systems far from equilibrium, in: 
{\it Phase Transitions and Critical Phenomena.} Vol. 19, 
C. Domb and J. Lebowitz (eds.), 
Academic Press, London (2001).

%
\bibitem{Schu20}
G.M. Sch\"utz, 
Duality from integrability: annihilating random walks with pair deposition,
J. Phys. A: Math. Theor. \textbf{53} 355003 (2020)


\bibitem{Simo09}
D. Simon,
Construction of a coordinate Bethe ansatz for the asymmetric simple exclusion process with open boundaries,
J. Stat. Mech. \textbf{2009}, P07017 (2009).

\bibitem{Spoh83}
H. Spohn,
Long-range correlations for stochastic lattice gases in a non-equilibrium
steady state.
J. Phys. A: Math. Gen. \textbf{16}, 4275--4291 (1983).

%
\bibitem{Sudb95}
A. Sudbury and P. Lloyd, 
Quantum operators in classical probability theory. II: The concept of 
duality in interacting particle systems, 
Ann. Probab. \textbf{23}(4), 1816--1830 (1995). 

\end{thebibliography}
\end{document}